\documentclass[3p,numbers,sort&compress]{elsarticle}

\usepackage{setspace}
\usepackage{amsmath}
\usepackage{amssymb}
\usepackage{amsthm}
\usepackage{bm}
\usepackage{graphicx}
\usepackage{caption}
\usepackage{subcaption}
\usepackage{siunitx}
\usepackage{float}
\usepackage{algorithm}
\usepackage{algpseudocode}
\usepackage{booktabs}
\usepackage{threeparttable}
\usepackage{romannum}
\usepackage{placeins}
\usepackage{microtype}
\usepackage[hidelinks]{hyperref}
\usepackage[nameinlink,capitalise,noabbrev]{cleveref}
\crefname{section}{Section}{Sections}
\Crefname{section}{Section}{Sections}
\crefname{subsection}{Section}{Sections}
\Crefname{subsection}{Section}{Sections}
\crefname{subsubsection}{Section}{Sections}
\Crefname{subsubsection}{Section}{Sections}
\crefname{chapter}{Chapter}{Chapters}
\Crefname{chapter}{Chapter}{Chapters}
\crefname{figure}{Fig.}{Figs.}
\Crefname{figure}{Fig.}{Figs.}
\crefname{table}{Table}{Tables}
\Crefname{table}{Table}{Tables}
\crefname{equation}{Eq.}{Eqs.}
\Crefname{equation}{Eq.}{Eqs.}
\newcommand{\figref}[1]{\cref{#1}}
\newcommand{\tabref}[1]{\cref{#1}}
\newcommand{\secref}[1]{\cref{#1}}

\newcommand{\eqnref}[1]{\cref{#1}}
\newcommand{\algoref}[1]{Algorithm~\ref{#1}}

\hypersetup{
  pdftitle={Distribution-Agnostic Robust Trajectory Optimization via Chance-Constrained Reinforcement Learning},
  pdfauthor={Yashdeep Chaudhary, Roberto Armellin, Harry Holt, Marco Sagliano},
  pdfkeywords={Robust control, Reinforcement Learning, Chance constraints, Spacecraft trajectory optimization, Covariance Steering, Interplanetary transfer, Powered descent},
  pdfsubject={math.OC},
  pdfdisplaydoctitle=true,
  pdfstartview={FitH},
  pdfcreator={LaTeX with hyperref}
}

\newcommand{\R}{\mathbb{R}}

\DeclareMathOperator{\E}{\mathbb{E}}

\DeclareMathOperator{\Cov}{Cov}
\DeclareMathOperator{\Prob}{\mathbb{P}}

\DeclareMathOperator{\diag}{diag}

\newcommand{\vect}[1]{\boldsymbol{#1}}
\newcommand{\mat}[1]{\boldsymbol{#1}}

\newcommand{\norm}[1]{\left\lVert #1 \right\rVert}
\newcommand{\abs}[1]{\left\lvert #1 \right\rvert}

\journal{arXiv}

\begin{document}

\begin{frontmatter}

    \title{Distribution-Agnostic Robust Trajectory Optimization via Chance-Constrained Reinforcement Learning}

    \author[auckland]{Yashdeep Chaudhary\corref{cor1}}
    \ead{ycha351@aucklanduni.ac.nz}

    \author[auckland]{Roberto Armellin}

    \author[esa]{Harry Holt}

    \author[unibo]{Marco Sagliano}

    \cortext[cor1]{Corresponding author}

    \affiliation[auckland]{
        organization={Faculty of Engineering, Waipapa Taumata Rau -- University of Auckland},
        addressline={20 Symonds Street},
        city={Auckland},
        postcode={1010},
        country={New Zealand}
    }

    \affiliation[esa]{
        organization={Advanced Concepts Team, European Space Research and Technology Centre, ESA},
        city={Noordwijk},
        postcode={2201 AZ},
        country={The Netherlands}
    }

    \affiliation[unibo]{
        organization={Department of Industrial Engineering, Alma Mater Studiorum - Università di Bologna},
        city={Via Zamboni 33, Bologna},
        postcode={40126},
        country={Italy}
    }

    \begin{abstract}
        This paper presents a distribution-agnostic robust trajectory-optimization framework based on chance-constrained reinforcement learning. The uncertainty is represented here through initial conditions and process noise, with the only requirement being that it can be sampled. A deterministic nominal trajectory is first computed offline, and reinforcement learning is then used only to robustify that baseline through a structured affine closed-loop correction law comprising a feedforward control adjustment and time-varying feedback gains. Probabilistic feasibility is enforced empirically through rollout-based upper-tail quantiles, while terminal dispersion is regulated through covariance-feasibility penalties. The framework is assessed on two materially different trajectory design problems. The flagship case study is a three-dimensional multi-impulse Earth--Mars transfer, where the learned policy is benchmarked against a recent robust trajectory-optimization reference under Gaussian uncertainty and then evaluated under bounded uniform uncertainty and under process disturbances not seen during training. The second case study is a stochastic atmospheric pinpoint rocket landing problem, used to assess portability to a short-horizon continuous-thrust setting with drag, mass depletion, and glide-slope constraints. The results show that the proposed framework can remain competitive in upper-tail fuel cost while preserving probabilistic feasibility, and that the same robustification scaffold can be carried across heterogeneous spacecraft trajectory planning problems without redesign of its core stochastic-control structure.
    \end{abstract}


\end{frontmatter}

\pagenumbering{arabic}

\section{Introduction}
\label{sec:introduction}

Spacecraft trajectory design is inherently shaped by uncertainty. Navigation errors, control-execution inaccuracies, unmodeled perturbations, atmospheric variability, and parametric mismatch can all drive a vehicle away from its nominal path, leading to fuel overruns, loss of terminal accuracy, or outright mission failure if not accounted for during the design phase~\cite{battin2001IntroductionMathematicsMethods,Conway2010-td,Betts1998-ey}. This issue becomes especially acute in missions with tight fuel margins and stringent terminal requirements, such as interplanetary rendezvous, pinpoint landing, and proximity operations, where the traditional workflow of computing a nominal trajectory and then appending conservative margins or late-stage feedback tuning is often either inefficient or insufficient.

Stochastic optimal control provides a principled framework for trajectory design under uncertainty by optimizing performance over a distribution of possible system evolutions rather than along a single nominal path~\cite{bertsekas1996StochasticOptimalControl}. Within this broader setting, chance-constrained formulations allow path and terminal requirements to be enforced probabilistically, thereby trading absolute conservatism for quantified risk~\cite{charnes1959ChanceConstrainedProgramming,prekopa1995StochasticProgramming}. Closely related covariance steering and covariance control methods further aim to regulate not only the mean trajectory but also the dispersion of the state distribution, which is particularly valuable when terminal accuracy is itself a mission requirement~\cite{Chen2016-go,ridderhof2021ChanceConstrainedCovarianceControl}. In the context of aerospace applications, these ideas have led to important advances in powered descent, low-thrust transfers, and robust interplanetary design~\cite{ridderhof2019MinimumfuelPoweredDescent,ridderhof2018UncertaintyQuanticationControl,benedikter2022ConvexApproachCovariance,benedikter2022ConvexApproachStochastic,marmo2025ChanceConstraintRobustTrajectory,Zhang2024-ih,Zhang2025-bj}.

These classical robust-optimization approaches are powerful, but their tractability often depends on structural assumptions that become restrictive in realistic mission settings. In particular, many formulations rely on Gaussian uncertainty models, linear or locally linear covariance propagation, and analytical reformulations of probabilistic constraints into deterministic surrogates. Such assumptions are often well-justified and effective within their intended regime, but they become increasingly difficult to maintain when uncertainty is bounded or non-Gaussian, when nonlinear dispersion growth is significant, or when constraint metrics are difficult to reformulate analytically. In those regimes, the resulting approximations may become overly conservative, brittle, or difficult to generalize across problem classes.

Reinforcement learning (RL) offers a complementary route. By learning policies directly from interaction with a stochastic simulation environment, RL can accommodate nonlinear dynamics, non-convex constraint geometries, and simulation-based uncertainty models without requiring closed-form reformulations of the underlying robust-control problem~\cite{sutton2018ReinforcementLearningIntroduction}. This flexibility has motivated growing interest in spacecraft applications including orbital transfers, entry, descent, and landing, proximity operations, and attitude control~\cite{zavoli2021ReinforcementLearningRobust,federici2021AutonomousGuidanceCislunar,gaudet2020DeepReinforcementLearning,Furfaro2020-pd,Federici2025-oc,Gaudet2020-cu,Hovell2021-jw,Federici2021-vg,Yuan2022-in,Mu2024-eo,Holt2025-uc,Capra2025-qe,vedant2019ReinforcementLearningSpacecraft,Zheng2021-md,Chen2023-em,chaudhary2024LowThrustCisLunarTransfers}. However, standard RL optimizes expected return rather than probabilistic feasibility, and unstructured end-to-end policies are often difficult to interpret, difficult to train efficiently, and weakly aligned with the structured optimization architectures typically preferred in safety-critical aerospace systems.

This paper develops a \emph{robust chance-constrained reinforcement learning} (RCCRL) framework for spacecraft trajectory optimization that combines a deterministic optimal-control backbone with sample-based probabilistic evaluation within RL. The central idea is to deliberately avoid replacing mature trajectory optimization machinery with an unconstrained learned controller. Instead, a deterministic nominal reference trajectory is first computed offline, in this work via Sequential Convex Programming (SCP), and the RL policy is then trained only to robustify that nominal solution under uncertainty. At each decision point, the policy outputs (i) a feedforward correction to the nominal control, and (ii) a time-varying linear feedback gain matrix that maps state deviation to a corrective action. Chance constraints are approximated empirically from rollout ensembles using upper-tail quantiles of mission-relevant metrics, while terminal dispersion is handled via a covariance feasibility requirement. This yields a closed-loop controller that remains structured and interpretable, yet can be trained using simulation-based uncertainty models beyond the Gaussian assumptions typically required for analytical tractability. In this paper, \emph{distribution-agnostic} refers to this control-design methodology, where the RCCRL framework does not require uncertainty to belong to a closed-form family admitting a bespoke analytical reformulation, provided it can be sampled in simulation.

Rather than claiming that RL-based approaches should supplant analytically grounded robust optimization whenever the latter is available, the intention of this paper is to demonstrate that the proposed RCCRL framework provides a structured simulation-based robustification architecture for regimes in which uncertainty handling through rollout ensembles and closed-loop policy adaptation are more natural than closed-form probabilistic reformulation. In this architecture, the nominal mission geometry is handled by a mature deterministic trajectory optimizer, whereas the learning component is restricted to robust closed-loop augmentation around that nominal solution. This preserves compatibility with established trajectory-design practice while avoiding the brittleness of unconstrained end-to-end policy learning.

The contributions of this paper are threefold. First, we formulate RCCRL as a finite-horizon stochastic trajectory-optimization framework in which probabilistic feasibility is enforced through sample-based quantile evaluation over Monte Carlo (MC) ensembles, thereby avoiding closed-form chance-constraint reformulations and permitting non-Gaussian uncertainty models so long as they can be sampled. Second, we adopt a structured affine policy parameterization that learns both feedforward control corrections and time-varying linear feedback gains around a deterministic nominal trajectory, allowing the controller to shape robustness explicitly while remaining interpretable. Third, building on the standard nominal--robust separation used in robust trajectory optimization, we formulate the RL component specifically as a closed-loop augmentation layer around a problem-appropriate reference trajectory. This design choice preserves the role of established trajectory optimization methods for generating the nominal mission geometry, while using RL to learn distribution-informed feedforward corrections and feedback gains under uncertainty.

To validate the proposed methodology across disparate problem settings, we benchmark the RCCRL framework on (i) a three-dimensional multi-impulse Earth--Mars heliocentric transfer, using the robust trajectory-optimization formulation of Marmo et al.~\cite{marmo2025ChanceConstraintRobustTrajectory} as a reference point under Gaussian uncertainty, and (ii) a stochastic atmospheric pinpoint rocket landing problem with continuous thrust, drag, and glide-slope constraints, informed by prior stochastic covariance-control and powered-descent formulations in Benedikter et al.~\cite{benedikter2022ConvexApproachStochastic} and Liu~\cite{liu2019FuelOptimalRocketLanding}. Among the two applications, the Earth--Mars transfer provides the most direct benchmark setting as it has the clearest connection to established robust astrodynamics methods and enables the proposed learning-based robustification strategy to be compared against a recent chance-constrained trajectory-optimization reference under comparable Gaussian uncertainty assumptions. The rocket-landing case study plays a complementary role by assessing whether the same RCCRL scaffold remains effective when carried into a markedly different dynamical regime involving short-horizon continuous-thrust control, aerodynamic effects, and state-dependent path constraints, thereby demonstrating that the framework is not tied to one specific transcription or dynamical regime.

The remainder of the paper is organized as follows. \secref{sec:framework} presents the RCCRL framework, including the stochastic problem setting, the structured affine control law, the empirical approximation of chance constraints, and the observation--action design used for policy learning. \secref{sec:cs1_em} instantiates the framework for the Earth--Mars transfer problem, and reports the corresponding results under Gaussian and bounded uniform uncertainty together with a robustness study under unmodeled process noise. \secref{sec:cs2_rktland} presents the atmospheric rocket-landing case study. \secref{sec:conclusion} concludes the paper with a synthesis of the methodological implications across both case studies.
\section{RCCRL Framework}
\label{sec:framework}

The RCCRL framework is built around three design choices. First, a nominal reference trajectory is generated offline using a problem-appropriate deterministic design method. Second, RL is used to augment that nominal design stage by learning a structured, affine, closed-loop policy around the reference. Third, probabilistic feasibility is handled empirically through MC ensemble rollouts and sample-based upper-tail statistics, thereby avoiding closed-form chance-constraint reformulations tied to a specific uncertainty family.

Throughout this section, the presentation is deliberately kept generic. The specific dynamical models, constraint metrics, normalization bounds, and training hyperparameters are introduced later for the Earth--Mars transfer and atmospheric rocket-landing case studies.

\subsection{Finite-horizon Stochastic Problem Setting}
\label{subsec:framework_problem}

Consider a discrete-time finite-horizon stochastic control problem over decision nodes $k = 0,\dots,N$, with state $\vect{x}_k \in \R^{n_x}$ and control $\vect{u}_k \in \R^{n_u}$. Uncertainty is represented through sampled initial conditions and, when applicable, sampled process disturbances. For the $i$-th sample in an MC ensemble of size $N_s$, the dynamics are written as
\begin{equation}
\vect{x}_{k+1}^{(i)} = f_k\!\left(\vect{x}_k^{(i)}, \vect{u}_k^{(i)}, \vect{w}_k^{(i)}\right),
\qquad
k = 0,\dots,N-1,\;\; i = 1,\dots,N_s,
\label{eq:rccrl_dynamics_generic}
\end{equation}
where $f_k(\cdot)$ denotes the problem-dependent discrete-time state transition and $\vect{w}_k^{(i)} \in \R^{n_w}$ is a disturbance realization.

The initial state ensemble is sampled from a prescribed distribution,
\begin{equation}
\vect{x}_0^{(i)} \sim \mathcal{D}_0,
\qquad i = 1,\dots,N_s,
\label{eq:rccrl_initial_distribution}
\end{equation}
and, when process uncertainty is included, disturbances are sampled according to
\begin{equation}
\vect{w}_k^{(i)} \sim \mathcal{D}_{w,k},
\qquad
k = 0,\dots,N-1,\;\; i = 1,\dots,N_s.
\label{eq:rccrl_process_distribution}
\end{equation}
The framework only assumes that the relevant uncertainty model can be sampled, which permits Gaussian, bounded uniform, and other simulation-defined uncertainty descriptions to be handled within the same training and evaluation scaffold. This is the sense in which RCCRL is termed \emph{distribution-agnostic}, where the control-design procedure does not rely on a closed-form uncertainty family or on an analytical chance-constraint reformulation tied to one.

At each decision node, the propagated ensemble provides an empirical approximation of the state distribution. Its first two moments (mean and covariance) are computed as
\begin{subequations}
\label{eq:rccrl_empirical_moments}
\begin{align}
\hat{\vect{x}}_k &= \frac{1}{N_s}\sum_{i=1}^{N_s}\vect{x}_k^{(i)},
\\
\hat{\mat{P}}_k &= \frac{1}{N_s-1}\sum_{i=1}^{N_s}
\left(\vect{x}_k^{(i)}-\hat{\vect{x}}_k\right)
\left(\vect{x}_k^{(i)}-\hat{\vect{x}}_k\right)^\top.
\end{align}
\end{subequations}
These empirical moments serve two purposes. First, they provide the compact uncertainty summary observed by the policy. Second, they supply the ensemble-level statistics needed for covariance-based terminal-dispersion assessment, without requiring an analytically tractable propagated distribution.

At a conceptual level, the stochastic problem may be written as
\begin{equation}
\begin{aligned}
\min_{\pi}\quad & \E\!\left[\sum_{k=0}^{N-1}\ell_k(\vect{x}_k,\vect{u}_k) + \ell_N(\vect{x}_N)\right] \\
\text{s.t.}\quad
& \vect{x}_{k+1} = f_k(\vect{x}_k,\vect{u}_k,\vect{w}_k), \\
& \vect{x}_0 \sim \mathcal{D}_0,\qquad \vect{w}_k \sim \mathcal{D}_{w,k}, \\
& \Prob\!\left(g_\ell\!\left(\{\vect{x}_k\},\{\vect{u}_k\}\right)\le 0\right)\ge 1-\beta_\ell,
\qquad \forall \ell = 1,\dots,N_c, \\
& h_m\!\left(\hat{\vect{x}}_N,\hat{\mat{P}}_N\right)\le 0,
\qquad \forall m = 1,\dots,N_e,
\end{aligned}
\label{eq:rccrl_stochastic_problem}
\end{equation}
where $\pi$ is the control policy, $\ell_k(\cdot)$ and $\ell_N(\cdot)$ are stage and terminal costs, $g_\ell(\cdot)$ denotes sample-wise chance-constrained quantities with risk levels $\beta_\ell \in [0,1]$, and $h_m(\cdot)$ denotes deterministic constraints imposed directly on ensemble-level terminal summaries such as covariance or dispersion. In practice, RCCRL does not attempt to solve \eqnref{eq:rccrl_stochastic_problem} through closed-form stochastic programming. Instead, it replaces the probabilistic problem with a sample-based penalized objective constructed from rollout ensembles. The point of the framework is therefore not to derive an analytical reformulation of the stochastic control problem, but to preserve its operational meaning while moving feasibility evaluation into simulation.

\subsection{Decoupled Nominal Backbone and Affine Robustification}
\label{subsec:framework_affine}

A defining feature of RCCRL is the decoupling between nominal trajectory design and robust closed-loop augmentation. The method assumes that a nominal reference trajectory and associated nominal controls,
\begin{equation}
\left\{\vect{x}_k^{\mathrm{nom}},\,\vect{u}_k^{\mathrm{nom}}\right\}_{k=0}^{N},
\label{eq:rccrl_nominal_reference}
\end{equation}
are available from a deterministic optimizer suited to the problem class. In this work, these references are generated via SCP, but the framework itself does not depend on SCP specifically. Any suitable nominal design procedure may be used, including other deterministic optimization pipelines or, if desired, even deterministic RL-based nominal design. The only requirement is the availability of a coherent reference trajectory around which robustification can be learned.

At each node $k$, the RCCRL policy outputs
\begin{equation}
\vect{a}_k =
\begin{bmatrix}
\Delta \vect{u}_k \\
\mathrm{vec}(\mat{K}_k)
\end{bmatrix},
\label{eq:rccrl_action_generic}
\end{equation}
where $\Delta \vect{u}_k \in \R^{n_u}$ is a feedforward correction to the nominal control, $\mat{K}_k \in \R^{n_u \times n_x}$ is a time-varying linear feedback gain matrix, and $\mathrm{vec}(\cdot)$ denotes the full vectorization operator.
The corrected nominal feedforward control is
\begin{equation}
\bar{\vect{u}}_k = \vect{u}_k^{\mathrm{nom}} + \Delta \vect{u}_k.
\label{eq:rccrl_corrected_nominal_control}
\end{equation}
A policy-adjusted reference trajectory is then propagated under disturbance-free dynamics:
\begin{equation}
\bar{\vect{x}}_{k+1} = f_k\!\left(\bar{\vect{x}}_k,\bar{\vect{u}}_k,\vect{0}\right),
\qquad
\bar{\vect{x}}_0 = \hat{\vect{x}}_0,
\label{eq:rccrl_reference_propagation}
\end{equation}
so that the feedforward correction and feedback law are defined relative to a reference consistent with the current learned control decision.

For each sample in the ensemble, the applied control is given by
\begin{equation}
\vect{u}_k^{(i)} =
\bar{\vect{u}}_k + \mat{K}_k\!\left(\vect{x}_k^{(i)}-\bar{\vect{x}}_k\right),
\qquad i = 1,\dots,N_s.
\label{eq:rccrl_affine_control}
\end{equation}
The RL action therefore combines a sample-independent feedforward correction $\Delta\vect{u}_k$ with a node-varying affine feedback gain $\mat{K}_k$. This affine structure is the main architectural restriction imposed on the learned controller. It is adopted here as a structured and interpretable feedback parameterization, consistent with the linear feedback laws commonly used in covariance-control and chance-constrained robust trajectory-optimization formulations~\cite{benedikter2022ConvexApproachStochastic,marmo2025ChanceConstraintRobustTrajectory}. The feedforward term allows the policy to introduce a local correction to the nominal control sequence when a systematic bias is beneficial, while the feedback term provides state-dependent contraction or reshaping of the uncertainty ensemble. In the present implementation, the feedforward component is therefore not intended to rediscover the nominal mission geometry from scratch, but to adjust the deterministic backbone within the neighborhood explored during robust training. Large departures from the reference solution remain discouraged through the fuel, terminal-feasibility, and covariance penalties, while the feedback term carries the primary role of regulating the ensemble dispersion. The learned controller consequently acts as a robustification layer around a nominal backbone rather than as an unconstrained end-to-end replacement for trajectory design.

\subsection{Empirical Approximation of Chance Constraints}
\label{subsec:framework_chance}

For each probabilistic constraint indexed by $\ell = 1,\dots,N_c$, let
\begin{equation}
z_\ell^{(i)} = \Psi_\ell\left(\{\vect{x}_k^{(i)}\}_{k=0}^{N},\{\vect{u}_k^{(i)}\}_{k=0}^{N-1}\right)
\label{eq:rccrl_constraint_metric_sample}
\end{equation}
denote a scalar sample-wise constraint metric constructed from the $i$-th trajectory realization. The corresponding chance constraint is satisfied when
\begin{equation}
\Prob(z_\ell \le 0) \ge 1-\beta_\ell.
\label{eq:rccrl_chance_constraint}
\end{equation}
Typical examples include path-wise control-effort limits, terminal miss-distance conditions, and other sample-wise mission-feasibility metrics.

Instead of relying on a closed-form analytical reformulation of \eqnref{eq:rccrl_chance_constraint} that may depend on specific distributional assumptions, RCCRL approximates the constraint empirically through upper-tail order statistics of the rollout ensemble. Let the values $\{z_\ell^{(i)}\}_{i=1}^{N_s}$ be sorted in non-decreasing order as
\begin{equation}
z_{\ell,(1)} \le z_{\ell,(2)} \le \cdots \le z_{\ell,(N_s)}.
\label{eq:rccrl_sorted_metrics}
\end{equation}
The empirical $(1-\beta_\ell)$-quantile is then defined as
\begin{equation}
\hat{q}_\ell(1-\beta_\ell) = z_{\ell,(j_\ell)},
\qquad
j_\ell = \left\lceil (1-\beta_\ell)N_s \right\rceil,
\label{eq:rccrl_empirical_quantile}
\end{equation}
and the chance constraint is approximated by the empirical condition
\begin{equation}
\hat{q}_\ell(1-\beta_\ell) \le 0.
\label{eq:rccrl_empirical_chance_constraint}
\end{equation}

This empirical quantile approximation is the central stochastic device in RCCRL. It allows the framework to accommodate non-Gaussian, bounded, or otherwise simulation-defined uncertainty without deriving a bespoke reformulation for each new uncertainty model or each new feasibility metric. The corresponding tradeoff is that the resulting guarantees are empirical rather than analytic. Their credibility therefore depends on the fidelity of the simulation environment, the representativeness of the assumed uncertainty model, and the accuracy of the finite-sample quantile estimates.

In addition to sample-wise chance constraints, some requirements are imposed directly on ensemble-level terminal statistics. In this paper, terminal dispersion is enforced through covariance-feasibility conditions of the form
\begin{equation}
c_m = h_m(\hat{\vect{x}}_N,\hat{\mat{P}}_N) \le 0,
\qquad m = 1,\dots,N_e,
\label{eq:rccrl_ensemble_constraints}
\end{equation}
where $h_m(\cdot)$ denotes a deterministic mapping from the terminal empirical mean and covariance to a scalar feasibility quantity. Chance constraints therefore act on upper-tail sample-wise behavior, whereas covariance feasibility acts directly on the terminal ensemble summary.

\subsection{Observation and Action Parameterization}
\label{subsec:framework_policy}

The policy is formulated as a finite-horizon Markov decision process in which the observation at node $k$ is a compact summary of the ensemble state together with nominal reference information. The generic observation vector is
\begin{equation}
\vect{o}_k =
\begin{bmatrix}
\hat{\vect{x}}_k \\
\mathrm{vec}_{\mathrm{up}}(\hat{\mat{P}}_k) \\
\vect{u}_k^{\mathrm{nom}} \\
\tau_k
\end{bmatrix},
\label{eq:rccrl_observation}
\end{equation}
where $\hat{\vect{x}}_k$ and $\hat{\mat{P}}_k$ are the empirical mean and covariance from \eqnref{eq:rccrl_empirical_moments}, $\mathrm{vec}_{\mathrm{up}}(\cdot)$ denotes the stacked upper-triangular entries of the symmetric covariance matrix, and $\tau_k \in [0,1]$ is the normalized time-to-go. Including $\vect{u}_k^{\mathrm{nom}}$ preserves awareness of the nominal control backbone, whereas the covariance summary supplies explicit information about the size and geometry of the current uncertainty cloud.

The policy does not observe individual MC samples directly. The ensemble is used internally by the environment to propagate uncertainty and to evaluate sample-based quantile penalties, but the policy itself acts only on the compact summary $\left\{\hat{\vect{x}}_k,\hat{\mat{P}}_k,\vect{u}_k^{\mathrm{nom}},\tau_k\right\}$. This keeps the input dimension fixed independently of ensemble size and ties the learned controller to quantities that could, in principle, be supplied by onboard estimation and nominal trajectory data. However, the present policy architecture remains moment-conditioned, i.e., distributions with identical first two moments but substantially different skewness, tails, or multimodal structure may not be distinguishable from the policy input alone. This limitation is partly mitigated by evaluating constraint satisfaction on the propagated sample ensemble over the evolved trajectory through empirical quantile penalties rather than inferring feasibility from covariance reduction alone.

The action is parameterized exactly as in \eqnref{eq:rccrl_action_generic}. In the present implementation, the policy outputs both the feedforward correction $\Delta\vect{u}_k$ and the entries of the feedback gain matrix $\mat{K}_k$. The corrected nominal feedforward control is then formed using \eqnref{eq:rccrl_corrected_nominal_control}, and the sample-wise control is applied through the affine law in \eqnref{eq:rccrl_affine_control}. The feedback gain matrix is fully parameterized, so the policy outputs all $n_u n_x$ entries of $\mat{K}_k$ in vectorized form. This provides sufficient expressive power to reshape or contract the ensemble in a state-dependent manner while preserving the affine control structure.

This observation--action structure also admits a plausible onboard interpretation. The MC ensemble is required during training and offline validation to expose the policy to uncertainty and to evaluate empirical chance-constraint penalties, but the deployed policy would not require access to all samples in flight. Instead, $\hat{\vect{x}}_k$ and $\hat{\mat{P}}_k$ could be supplied by an onboard navigation or state-estimation filter, while $\vect{u}_k^{\mathrm{nom}}$ and the nominal reference trajectory are precomputed offline. At execution time, the policy would output $\Delta\vect{u}_k$ and $\mat{K}_k$ from the compact observation, construct the policy-adjusted feedforward command, and apply the affine feedback correction using the current estimated state deviation relative to the policy-adjusted reference state $\bar{\vect{x}}_k$.

The policy $\pi_\theta(\vect{a}_k \mid \vect{o}_k)$ is implemented as a stochastic actor--critic model trained with Proximal Policy Optimization (PPO). PPO is adopted here for practical rather than ideological reasons: the action space is continuous, the optimization is episodic and penalty-based, and stable policy-gradient updates are needed under repeated rollout sampling. The actor outputs the mean of a diagonal Gaussian action distribution, together with learned log-standard deviations, while the critic estimates the state value associated with the same observation input. The exact network widths, normalization ranges, and PPO hyperparameters are case-study-dependent and are therefore specified later with the Earth--Mars and atmospheric rocket-landing instantiations.

Consistent with the observation that implementation details can materially affect reported reinforcement-learning performance~\cite{engstrom2020ImplementationMattersDeep}, the RCCRL policies in this paper are trained using established and well-tested software libraries. Specifically, \texttt{Stable-Baselines3}~\cite{stable-baselines3} provides the PPO implementation used here, built on \texttt{PyTorch}~\cite{paszke2019PyTorchImperativeStyle}, while the training environments are implemented as custom \texttt{Gymnasium}~\cite{towers2024gymnasium}-compatible classes. This pipeline supports reproducibility, enables efficient vectorized rollout collection, and helps isolate algorithmic and environment-design effects from implementation-specific artifacts.

\subsection{Penalized RL Objective and Training Workflow}
\label{subsec:framework_objective}

For a rollout generated under policy $\pi_\theta$, let $\{M_j\}_{j=1}^{N_m}$ denote a set of scalar performance metrics extracted from the propagated ensemble trajectories and controls. These metrics capture mission objectives such as fuel expenditure, terminal accuracy, or dispersion-related behavior. An aggregate performance cost is then written as
\begin{equation}
J_{\mathrm{perf}} = \sum_{j=1}^{N_m}\alpha_j M_j,
\label{eq:rccrl_performance_cost}
\end{equation}
where $\alpha_j \ge 0$ are user-chosen weighting coefficients.

Chance-constraint violations are penalized through empirical quantiles of sample-wise constraint residuals:
\begin{equation}
\Phi_{\mathrm{cc}}
=
\sum_{\ell=1}^{N_c}
\lambda_\ell
\left[\hat{q}_\ell(1-\beta_\ell)\right]_{+}^{p_\ell},
\label{eq:rccrl_quantile_penalty}
\end{equation}
where $[\cdot]_+ = \max(\cdot,0)$, $\lambda_\ell > 0$ is a penalty weight, and $p_\ell \ge 1$ is a shaping exponent. This term can be used for any constraint that is naturally evaluated on a per-sample basis, including per-impulse control limits, path constraints, and terminal miss-distance or capture constraints. Thus, terminal feasibility need not be required for every propagated sample, and when appropriate, can be imposed probabilistically by including the relevant terminal residual in $\Phi_{\mathrm{cc}}$.
By contrast, deterministic ensemble-level feasibility conditions are penalized separately as
\begin{equation}
\Phi_{\mathrm{ens}}
=
\sum_{m=1}^{N_e}
\eta_m \left[c_m\right]_{+}^{r_m},
\label{eq:rccrl_ensemble_penalty}
\end{equation}
with $\eta_m > 0$ and $r_m \ge 1$. These terms are used for constraints that are defined directly on ensemble statistics rather than on individual samples, such as terminal covariance bounds, covariance-shape requirements, or deterministic terminal mean errors.

The total per-episode penalized cost is therefore
\begin{equation}
J_{\mathrm{tot}}
=
J_{\mathrm{perf}}
+
\Phi_{\mathrm{cc}}
+
\Phi_{\mathrm{ens}},
\label{eq:rccrl_total_cost}
\end{equation}
and the scalar episodic RL return supplied to PPO is taken as its negative:
\begin{equation}
R = -J_{\mathrm{tot}}.
\label{eq:rccrl_return}
\end{equation}
For finite-horizon problems, this return may be assigned terminally or decomposed into stepwise rewards, whose sum recovers the same episodic quantity. In either case, PPO maximizes the expected return induced by the sampled initial conditions, sampled disturbances, and policy stochasticity.

\eqnref{eq:rccrl_return} makes explicit that RCCRL does not enforce the original stochastic constraints as hard constraints within RL. Instead, reward design acts as the scalarization through which nominal performance, upper-tail probabilistic feasibility, and terminal covariance quality are negotiated during policy learning. This is both a strength and a limitation as it provides flexibility across problem classes, but it also shifts part of the method design burden onto the choice of metrics, penalty weights, and normalization.

\begin{algorithm*}[htbp]
\caption{RCCRL Training Workflow}
\label{alg:rccrl_workflow}
\begin{algorithmic}[1]
\Require Nominal design method, uncertainty models $\mathcal{D}_0$ and $\mathcal{D}_{w,k}$, ensemble size $N_s$, policy $\pi_\theta$
\Ensure Trained policy $\pi_\theta$ and nominal reference $\{\vect{x}_k^{\mathrm{nom}}, \vect{u}_k^{\mathrm{nom}}\}_{k=0}^{N}$

\State Compute the nominal reference trajectory $\{\vect{x}_k^{\mathrm{nom}}, \vect{u}_k^{\mathrm{nom}}\}_{k=0}^{N}$ using a problem-appropriate nominal design method
\State Initialize the policy and value-function parameters
\For{each training iteration}
    \State Sample an initial ensemble from $\mathcal{D}_0$
    \State Sample process disturbances from $\mathcal{D}_{w,k}$, if applicable
    \For{$k = 0,\dots,N-1$}
        \State Compute $\hat{\vect{x}}_k$ and $\hat{\mat{P}}_k$ from the propagated ensemble
        \State Form the observation $\vect{o}_k$ using \eqnref{eq:rccrl_observation}
        \State Sample the action $\vect{a}_k = [\Delta \vect{u}_k^\top,\mathrm{vec}(\mat{K}_k)^\top]^\top \sim \pi_\theta(\cdot \mid \vect{o}_k)$
        \State Compute $\bar{\vect{u}}_k$ and propagate $\bar{\vect{x}}_{k+1}$
        \State Apply the affine law in \eqnref{eq:rccrl_affine_control} to each sample
        \State Propagate the ensemble to node $k+1$
    \EndFor
    \State Evaluate the performance metrics, empirical quantiles, and covariance-feasibility penalties
    \State Compute the episodic return using \eqnref{eq:rccrl_return}
    \State Update the policy parameters using PPO
\EndFor
\end{algorithmic}
\end{algorithm*}

The decoupled RCCRL workflow is summarized in \algoref{alg:rccrl_workflow}. A deterministic reference trajectory is first computed offline using a problem-appropriate nominal design method. That reference is then embedded into the stochastic environment used for policy learning, where ensemble propagation, empirical quantile evaluation, and PPO-based policy updates are carried out iteratively.

The presented framework is instantiated next in two materially different optimization problems, with only the problem-specific dynamics, metrics, and scaling choices altered: first as a benchmark Earth--Mars transfer case study against an established robust trajectory-optimization reference, and then as an atmospheric pinpoint rocket-landing case study used to assess methodological portability beyond the interplanetary impulsive-transfer regime.
\section{Case Study~\Romannum{1}: Earth--Mars Transfer}
\label{sec:cs1_em}

The RCCRL framework is instantiated on a three-dimensional, time-fixed Earth--Mars heliocentric transfer under uncertainty. This mission class provides a long-horizon impulsive trajectory design problem with meaningful terminal-capture requirements and a direct point of comparison against the recent robust trajectory-optimization formulation of Marmo et al.~\cite{marmo2025ChanceConstraintRobustTrajectory}.
In addition to the Gaussian initial-uncertainty model, the same formulation is also evaluated under bounded uniform uncertainty with matched first and second moments, thereby probing the behavior of the learned affine correction law beyond the Gaussian regime.

\subsection{Problem Setup}
\label{subsec:em_dynamics}

The spacecraft is modeled as a point mass evolving under heliocentric two-body dynamics in the ecliptic J2000 frame. The state is
\begin{equation}
\vect{x}(t)
=
\begin{bmatrix}
\vect{r}(t) \\
\vect{v}(t)
\end{bmatrix}
\in \R^{6},
\label{eq:em_state}
\end{equation}
where $\vect{r}(t),\vect{v}(t)\in\R^{3}$ denote inertial position and velocity, respectively. Between maneuver nodes, the dynamics are purely Keplerian:
\begin{equation}
\dot{\vect{r}} = \vect{v},
\qquad
\dot{\vect{v}} = -\mu_\odot \frac{\vect{r}}{\norm{\vect{r}}^{3}},
\label{eq:em_twobody}
\end{equation}
with $\mu_\odot$ denoting the solar gravitational parameter.

The transfer departs Earth with zero hyperbolic excess velocity at epoch $t_0$ and must rendezvous with Mars, again with zero excess velocity, at the prescribed arrival time $t_f$. The deterministic boundary conditions are
\begin{equation}
\vect{x}^{-}(t_0) = \vect{x}_E(t_0),
\qquad
\vect{x}^{+}(t_f) = \vect{x}_M(t_f),
\label{eq:em_boundary_conditions}
\end{equation}
where $\vect{x}_E(t)$ and $\vect{x}_M(t)$ are the heliocentric ephemeris states of Earth and Mars, obtained here from the JPL DE405 ephemerides~\cite{Standish1998}, and the superscripts $-$ and $+$ denote the lower and upper limits of the state at the corresponding epoch.

\begin{table}[h!tb]
\centering
\caption{Problem parameters for the Earth--Mars transfer case study.}
\label{tab:em_problem_data}
\begin{threeparttable}
\footnotesize
\begin{tabular}{@{}lcc@{}}
\toprule
Parameter & Value & Unit \\
\midrule
\multicolumn{3}{@{}l}{\textbf{Normalization reference}} \\
\quad Reference length scale, $r_{\mathrm{scale}}$ & $1.4960\times10^{8}$ & km \\
\quad Solar gravitational parameter, $\mu_\odot$ & $1.3271\times10^{11}$ & km$^{3}$/s$^{2}$ \\
\quad Reference velocity scale, $v_{\mathrm{scale}}$ & $29.7845$ & km/s \\
\quad Reference time scale, $t_{\mathrm{scale}}$ & $5.023\times10^{6}$ & s \\
\midrule
\multicolumn{3}{@{}l}{\textbf{Discretization and transfer setup}} \\
\quad Total time of flight, $\Delta T$ & $348.7900$ & days \\
\quad Number of ballistic segments, $N$ & $20$ & -- \\
\quad Segment duration, $\Delta t$ & $17.4395$ & days \\
\quad Maximum per-impulse magnitude, $\Delta v_{\max}$ & $0.7600$ & km/s \\
\quad Risk level, $\beta$ & $5\%$ & -- \\
\midrule
\multicolumn{3}{@{}l}{\textbf{Boundary states}} \\
\quad Initial position, $\vect{r}^{-}(t_0)$ &
$[-140{,}699{,}693,\,-51{,}614{,}428,\,980]^\top$ & km \\
\quad Initial velocity, $\vect{v}^{-}(t_0)$ &
$[9.7746,\,-28.0783,\,4.3377\times10^{-4}]^\top$ & km/s \\
\quad Target position, $\vect{r}^{+}(t_f)$ &
$[-172{,}682{,}023,\,176{,}959{,}469,\,7{,}948{,}912]^\top$ & km \\
\quad Target velocity, $\vect{v}^{+}(t_f)$ &
$[-16.4274,\,-14.8605,\,9.2149\times10^{-2}]^\top$ & km/s \\
\midrule
\multicolumn{3}{@{}l}{\textbf{Initial and target uncertainty}} \\
\quad Initial position std. dev., $\sigma_r^{-}(t_0)$ & $1.5000\times10^{6}$ & km \\
\quad Initial velocity std. dev., $\sigma_v^{-}(t_0)$ & $9.4128\times10^{-2}$ & km/s \\
\quad Target position std. dev., $\sigma_r^{+}(t_f)$ & $1.5000\times10^{5}$ & km \\
\quad Target velocity std. dev., $\sigma_v^{+}(t_f)$ & $9.4128\times10^{-3}$ & km/s \\
\midrule
\multicolumn{3}{@{}l}{\textbf{Process-noise intensities}} \\
\quad Position noise intensity, $\sigma_{w,r}$ & $3.3466\times10^{5}$ & km \\
\quad Velocity noise intensity, $\sigma_{w,v}$ & $2.1048\times10^{-2}$ & km/s \\
\bottomrule
\end{tabular}
\end{threeparttable}
\end{table}

The trajectory is transcribed as a sequence of $N$ ballistic arcs of equal duration separated by impulsive maneuvers at the grid points
\begin{equation}
t_0 \le t_1 \le \cdots \le t_N = t_f,
\qquad
\Delta t = \frac{t_f - t_0}{N}.
\label{eq:em_time_grid}
\end{equation}
Let $\Delta \vect{v}_k \in \R^{3}$ denote the impulsive velocity change applied at node $k$. The jump conditions are
\begin{equation}
\vect{r}^{+}(t_k) = \vect{r}^{-}(t_k),
\qquad
\vect{v}^{+}(t_k) = \vect{v}^{-}(t_k) + \Delta \vect{v}_k,
\qquad
k = 0,\dots,N,
\label{eq:em_jump_conditions}
\end{equation}
and each impulse is bounded as
\begin{equation}
\norm{\Delta \vect{v}_k} \le \Delta v_{\max},
\qquad
k = 0,\dots,N.
\label{eq:em_dv_cap}
\end{equation}
The corresponding deterministic impulse cost is
\begin{equation}
J_{\Delta v} = \sum_{k=0}^{N} \norm{\Delta \vect{v}_k}.
\label{eq:em_dv_cost}
\end{equation}

All computations are performed in non-dimensional units to improve numerical conditioning:
\begin{equation}
\tilde{\vect{r}} = \frac{\vect{r}}{r_{\mathrm{scale}}},
\qquad
\tilde{\vect{v}} = \frac{\vect{v}}{v_{\mathrm{scale}}},
\qquad
\tilde{t} = \frac{t}{t_{\mathrm{scale}}},
\label{eq:em_nondim}
\end{equation}
with 1 AU used as the reference length and the associated circular heliocentric speed as the reference velocity, so that the non-dimensional solar parameter is unity by construction.
The values reported in \tabref{tab:em_problem_data} follow the benchmark configuration used in the Earth--Mars robust trajectory-optimization literature, which enables a direct Gaussian-case comparison in the results section.\footnote{The discrepancy between the adopted process-noise intensities and the printed values in Ref.~\cite{marmo2025ChanceConstraintRobustTrajectory} was confirmed in personal communication with the authors (November 2025), and the values adopted here are those consistent with the benchmark numerical results reported therein.}

\subsection{Uncertainty Model and Probabilistic Requirements}
\label{subsec:em_uncertainty}

At departure, the uncertain initial state is modeled as
\begin{equation}
\vect{x}(t_0) \sim \mathcal{D}_0,
\qquad
\E[\vect{x}(t_0)] = \vect{x}_0,
\qquad
\Cov[\vect{x}(t_0)] = \mat{P}_0,
\label{eq:em_initial_distribution}
\end{equation}
with nominal departure state
\begin{equation}
\vect{x}_0 = \vect{x}_E(t_0),
\label{eq:em_nominal_departure}
\end{equation}
and covariance
\begin{equation}
\mat{P}_0
=
\diag\!\left(
\sigma_r^{-}(t_0)^2 \mat{I}_3,\,
\sigma_v^{-}(t_0)^2 \mat{I}_3
\right).
\label{eq:em_initial_cov}
\end{equation}

Two initial-state models are considered.
First, in the \emph{Gaussian case},
\begin{equation}
\mathcal{D}_0 = \mathcal{N}(\vect{x}_0,\mat{P}_0),
\label{eq:em_initial_gaussian}
\end{equation}
which matches the uncertainty formulation typically adopted in existing chance-constrained and covariance-control studies, including~\cite{marmo2025ChanceConstraintRobustTrajectory}.
Second, in the \emph{uniform case}, each state component is sampled independently from a bounded uniform distribution with matched mean and marginal variance:
\begin{equation}
x_i(t_0)
\sim
\mathcal{U}\!\left(\mu_i-\sqrt{3}\,\sigma_i,\;\mu_i+\sqrt{3}\,\sigma_i\right),
\label{eq:em_initial_uniform}
\end{equation}
where $\mu_i$ and $\sigma_i$ are the corresponding entries of $\vect{x}_0$ and $\sqrt{\diag(\mat{P}_0)}$. This preserves the first two moments while changing the support and higher-order structure of the uncertainty.

In addition to initial-state dispersion, the framework also supports continuous-time process noise acting on the full state:
\begin{equation}
d\vect{x}(t) = f(\vect{x}(t))\,dt + \mat{G}\,d\vect{w}(t),
\label{eq:em_sde}
\end{equation}
with standard Wiener process $\vect{w}(t)$ and diffusion matrix
\begin{equation}
\mat{G}
=
\diag\!\left(
\sigma_{w,r}\mat{I}_3,\,
\sigma_{w,v}\mat{I}_3
\right).
\label{eq:em_diffusion}
\end{equation}
Within the RCCRL environment, this disturbance is discretized segment-wise using the local state-transition matrix, yielding an additive process-noise covariance applied to the propagated ensemble. In the main training campaign, process noise is disabled during training so that the policy is learned primarily against initial-state uncertainty. It is then reintroduced at evaluation time to test robustness under disturbance mismatch. This is intentional and is studied explicitly in \secref{subsec:em_results_noise}.

The Earth--Mars transfer is subject to three probabilistic requirements. First, the impulse magnitude must satisfy
\begin{equation}
\Prob\!\left(\norm{\Delta \vect{v}_k} \le \Delta v_{\max}\right) \ge 1-\beta,
\qquad
k = 0,\dots,N.
\label{eq:em_prob_dv}
\end{equation}
Second, the terminal position must lie within the Mars sphere of influence (SOI), of radius $r_{\mathrm{SOI}}$:
\begin{equation}
\Prob\!\left(e_r \le r_{\mathrm{SOI}}\right) \ge 1-\beta,
\qquad
e_r = \norm{\vect{r}(t_f)-\vect{r}_M(t_f)}.
\label{eq:em_prob_capture}
\end{equation}
Third, the terminal covariance must remain below a prescribed target covariance $\mat{P}(t_f) \preceq \mat{P}_f,$ with
\begin{equation}
\mat{P}_f
=
\diag\!\left(
\sigma_r^{+}(t_f)^2 \mat{I}_3,\,
\sigma_v^{+}(t_f)^2 \mat{I}_3
\right).
\label{eq:em_target_cov}
\end{equation}

A common risk level $\beta$ is used for all empirical chance constraints in this case study to preserve a compact benchmark comparison. The RCCRL framework itself does not require a uniform risk allocation, and separate risk levels may be assigned to individual constraints when mission requirements warrant different violation probabilities. In the learning environment, the impulse cap and SOI-capture requirement are handled through empirical upper-tail penalties as described in \secref{subsec:framework_chance}, while the terminal covariance constraint is assessed through the symmetric difference $\Delta \mat{P} = \mat{P}_f - \hat{\mat{P}}_N$ and the scalar violation measure
\begin{equation}
\varepsilon_{\mathrm{cov}}
=
\sum_{\Lambda_j < 0}
\abs{\Lambda_j(\Delta \mat{P})},
\label{eq:em_cov_violation}
\end{equation}
where $\Lambda_j(\cdot)$ denotes the $j$-th eigenvalue. This scalar penalty is zero if and only if the terminal covariance requirement is satisfied.

\subsection{Deterministic Nominal Baseline}
\label{subsec:em_scp}

Before introducing the learning-based robustness layer, a deterministic fuel-optimal reference trajectory is computed offline using SCP. This baseline provides the nominal state--control sequence around which the affine RCCRL closed-loop augmentation is learned and later applied.

The deterministic problem minimizes the total impulse cost \eqnref{eq:em_dv_cost} subject to the two-body dynamics, the boundary conditions \eqnref{eq:em_boundary_conditions}, and the per-node impulse cap \eqnref{eq:em_dv_cap}. In non-dimensional form, the resulting nominal optimal-control problem is
\begin{subequations}
\label{eq:em_scp_problem}
\begin{align}
\min_{\{\Delta \vect{v}_k\}_{k=0}^{N}}
\quad &
\sum_{k=0}^{N} \norm{\Delta \vect{v}_k}
\label{eq:em_scp_obj}
\\
\text{s.t.}\quad &
\tilde{\vect{x}}_0^{-} = \tilde{\vect{x}}_E(t_0),
\qquad
\tilde{\vect{x}}_{N}^{+} = \tilde{\vect{x}}_M(t_f),
\label{eq:em_scp_bc}
\\
&
\tilde{\vect{x}}_{k+1}^{-}
=
\tilde{\mathcal{F}}_k\!\left(\tilde{\vect{x}}_k^{+},\Delta \tilde{t}\right),
\qquad
k = 0,\dots,N-1,
\label{eq:em_scp_dyn}
\\
&
\tilde{\vect{x}}_k^{+}
=
\begin{bmatrix}
\tilde{\vect{r}}_k^{-} \\
\tilde{\vect{v}}_k^{-} + \Delta \tilde{\vect{v}}_k
\end{bmatrix},
\qquad
k = 0,\dots,N,
\label{eq:em_scp_jump}
\\
&
\norm{\Delta \vect{v}_k} \le \Delta v_{\max},
\qquad
k = 0,\dots,N.
\label{eq:em_scp_cap}
\end{align}
\end{subequations}
Here, $\tilde{\mathcal{F}}_k(\cdot)$ denotes the normalized ballistic propagation operator over one segment. 

The nonlinear dynamics are successively linearized about a reference impulsive trajectory, and the required discrete-time Jacobians are obtained by automatic differentiation through the Runge--Kutta propagation scheme. Each convex subproblem is posed in \texttt{CVXPY}~\cite{diamond2016cvxpy,agrawal2018rewriting} and solved using \texttt{MOSEK}~\cite{mosek} until the defect between the affine model and full nonlinear propagation falls below the prescribed tolerance. The initial reference is generated from a Lambert transfer between the Earth and Mars ephemeris positions at $(t_0,t_f)$, with departure and arrival impulses initialized relative to the planetary ephemeris velocities.

\begin{table}[h!tb]
\centering
\caption{Deterministic SCP baseline results for the Earth--Mars transfer case study.}
\label{tab:em_scp_results}
\begin{threeparttable}
\footnotesize
\begin{tabular}{@{}lcc@{}}
\toprule
Metric & Value & Unit \\
\midrule
Total deterministic $\Delta v$ cost & $10.0585$ & km/s \\
Terminal position error (nonlinear validation) & $2.3240$ & km \\
Terminal velocity error (nonlinear validation) & $1.6376\times10^{-7}$ & km/s \\
Number of SCP iterations & $6$ & -- \\
\bottomrule
\end{tabular}
\end{threeparttable}
\end{table}

\pagebreak

The resulting deterministic baseline, summarized in \tabref{tab:em_scp_results}, is accurate to a level that is negligible relative to both the Mars SOI radius and the target terminal-dispersion scales. It therefore provides a suitable nominal anchor for robust policy learning.

\subsection{RCCRL Instantiation}
\label{subsec:em_rccrl_instantiation}

\subsubsection{Control Parameterization and Observation Design}
\label{subsubsec:em_policy_design}

For the Earth--Mars case study, the nominal control backbone is the SCP impulse sequence $\left\{\Delta \vect{v}_k^{\mathrm{ref}}\right\}_{k=0}^{N-1}.$
At each node, the RCCRL policy learns a feedforward corrective impulse $\Delta \vect{v}_k^{\mathrm{corr}}\in\R^{3}$ and a time-varying feedback gain matrix $\mat{K}_k\in\R^{3\times 6}$. For the $i$-th ensemble member, the applied impulse is
\begin{equation}
\Delta \vect{v}_k^{(i)}
=
\Delta \vect{v}_k^{\mathrm{ref}}
+
\Delta \vect{v}_k^{\mathrm{corr}}
+
\mat{K}_k\!\left(\vect{x}_k^{(i)}-\bar{\vect{x}}_k\right),
\label{eq:em_affine_impulse}
\end{equation}
where $\bar{\vect{x}}_k$ is the policy-adjusted reference state from \secref{subsec:framework_affine}. The actor therefore operates on a $21$-dimensional action space: $3$ feedforward-control entries and $18$ gain entries.

The observation combines the current ensemble summary, the nominal control reference, and the remaining horizon information:
\begin{equation}
    \vect{o}_k = \left[\hat{\vect{x}}_k,\,\mathrm{vec}_{\mathrm{up}}(\hat{\mat{P}}_k),\,\Delta \vect{v}_k^{\mathrm{ref}},\,\tau_k\right]^\top \in \R^{31},        
\end{equation}
where the $31$ entries consist of $6$ mean-state components, $21$ upper-triangular covariance terms, $3$ reference-control entries, and $1$ normalized time-to-go term. Each observation component is linearly scaled to $[-1,1]$ using fixed bounds obtained from precomputed MC sweeps around the nominal trajectory.

In the terminal segment, a bi-impulsive construction is used. The first leg is supplied by the policy through \eqnref{eq:em_affine_impulse}, after which a deterministic second-leg correction $\Delta \vect{v}_{N}^{\mathrm{leg2}}$ is applied uniformly across the ensemble so that the nominal terminal velocity matches the required Mars-arrival state. This second-leg correction contributes to both the total $\Delta v$ cost and the terminal penalties.

\subsubsection{Reward Construction}
\label{subsubsec:em_reward}

The Earth--Mars reward instantiates the generic penalized RCCRL objective using percentile-based fuel and feasibility terms. At each non-terminal node, the frequent (per-step) reward is
\begin{subequations}
\label{eq:em_frequent_reward}
\begin{align}
r_k^{\mathrm{freq}}
&=
r_k^{\Delta v}+r_k^{\mathrm{viol}},
\label{eq:em_freq_total}
\\
r_k^{\Delta v}
&=
-\lambda_{\Delta v}\,
\hat{q}_{p}\!\left(\left\{\norm{\Delta \vect{v}_k^{(i)}}\right\}_{i=1}^{N_s}\right),
\label{eq:em_freq_dv}
\\
r_k^{\mathrm{viol}}
&=
-\lambda_{\mathrm{viol}}\,
\max\!\left\{
0,\,
\hat{q}_{p}\!\left(\left\{\norm{\Delta \vect{v}_k^{(i)}}-\Delta v_{\max}\right\}_{i=1}^{N_s}\right)
\right\}.
\label{eq:em_freq_viol}
\end{align}
\end{subequations}

At the terminal node, additional penalties are introduced for the deterministic second-leg correction, SOI-capture violation, and covariance infeasibility:
\begin{subequations}
\label{eq:em_terminal_reward}
\begin{align}
r_N^{\mathrm{term}}
&=
r_N^{\Delta v,\mathrm{leg2}}
+r_N^{\mathrm{viol,leg2}}
+r_N^{\mathrm{pos}}
+r_N^{\mathrm{cov}}
+r_N^{\mathrm{bonus}},
\label{eq:em_term_total}
\\
r_N^{\Delta v,\mathrm{leg2}}
&=
-\lambda_{\Delta v}\norm{\Delta \vect{v}_N^{\mathrm{leg2}}},
\label{eq:em_term_leg2_cost}
\\
r_N^{\mathrm{viol,leg2}}
&=
-\lambda_{\mathrm{viol}}
\max\!\left\{0,\norm{\Delta \vect{v}_N^{\mathrm{leg2}}}-\Delta v_{\max}\right\},
\label{eq:em_term_leg2_viol}
\\
r_N^{\mathrm{pos}}
&=
-\lambda_{\mathrm{pos}}
\min\!\left\{
\frac{
\max\!\left\{
0,\,
\hat{q}_{p}\!\left(\left\{\norm{\vect{r}_N^{(i)}-\vect{r}_M(t_f)}\right\}_{i=1}^{N_s}\right)-r_{\mathrm{SOI}}
\right\}
}{r_{\mathrm{SOI}}},
\,
c_{\mathrm{pos}}
\right\},
\label{eq:em_term_pos}
\\
r_N^{\mathrm{cov}}
&=
-\lambda_{\mathrm{cov}}\,
\varepsilon_{\mathrm{cov}}.
\label{eq:em_term_cov}
\end{align}
\end{subequations}

A small positive bonus is awarded when the terminal metrics simultaneously satisfy tight tolerances, encouraging the policy to enter the feasible regime before pursuing marginal fuel improvements. The full reward configuration is listed in \tabref{tab:em_reward_config}.

\begin{table}[htbp]
\centering
\caption{Reward parameters for the Earth--Mars transfer case study.}
\label{tab:em_reward_config}
\begin{threeparttable}
\footnotesize
\begin{tabular}{@{}lcc@{}}
\toprule
Parameter & Value & Unit \\
\midrule
\multicolumn{3}{@{}l}{\textbf{Reward weights}} \\
\quad $\Delta v$ cost weight, $\lambda_{\Delta v}$ & $40.0$ & -- \\
\quad $\Delta v$ violation weight, $\lambda_{\mathrm{viol}}$ & $400.0$ & -- \\
\quad Terminal position-dispersion weight, $\lambda_{\mathrm{pos}}$ & $100.0$ & -- \\
\quad Terminal covariance weight, $\lambda_{\mathrm{cov}}$ & $5.0\times10^{7}$ & -- \\
\quad Terminal bonus weight, $\lambda_{\mathrm{bonus}}$ & $180.0$ & -- \\
\midrule
\multicolumn{3}{@{}l}{\textbf{Terminal tolerances}} \\
\quad $\Delta v$ violation tolerance & $1.0\times10^{-2}$ & km/s \\
\quad Position-dispersion tolerance & $2.885\times10^{4}$ & km \\
\quad Covariance-violation tolerance & $1.0\times10^{-6}$ & -- \\
\quad Position penalty cap, $c_{\mathrm{pos}}$ & $500.0$ & -- \\
\bottomrule
\end{tabular}
\end{threeparttable}
\end{table}

The $\Delta v$- and terminal-position-related terms are evaluated in dimensional units after rescaling the non-dimensional states and controls, whereas covariance-related terms are evaluated in normalized form. This mixed treatment preserves physical interpretability for fuel and capture penalties while maintaining numerical stability for the covariance-feasibility term.

\subsubsection{Policy Architecture and Training Configuration}
\label{subsubsec:em_ppo_setup}

The Earth--Mars RCCRL policy is trained with PPO using the generic actor--critic structure introduced in \secref{subsec:framework_policy}. For this case study, the actor contains approximately $40{,}500$ trainable parameters and the critic approximately $6{,}800$, for a total of about $47{,}300$ parameters. The network widths and PPO hyperparameters were selected to balance exploration, sample efficiency, and training stability in the long-horizon sparse-reward regime induced by the impulsive transfer problem.

The full training configuration is reported in \tabref{tab:em_ppo_config}. The discount factor is set very close to unity to account for the long horizon, the learning rate and clip range are linearly annealed during training, and the training ensemble size is $N_s^{\mathrm{train}}=512$. The reported campaign uses $8$ parallel environments and a total of $8.0\times10^{7}$ training timesteps.
The value $N_s^{\mathrm{train}}=512$ was chosen as a compromise between tail-statistic fidelity and rollout cost. This ensemble size provides a sufficiently populated upper tail for estimating the empirical 95th-percentile penalties during training, while still allowing repeated PPO updates with $8$ parallel environments. The final reported statistics are recomputed using $N_s^{\mathrm{eval}}=100{,}000$ independent MC samples.

\begin{table}[h!tb]
\centering
\caption{PPO architecture and training parameters for the Earth--Mars transfer case study.}
\label{tab:em_ppo_config}
\begin{threeparttable}
\footnotesize
\begin{tabular}{@{}lcc@{}}
\toprule
Parameter & Symbol & Value \\
\midrule
\multicolumn{3}{@{}l}{\textbf{Network architecture}} \\
\quad Observation length & $n_o$ & $31$ \\
\quad Action length & $n_a$ & $21$ \\
\quad Number of hidden layers & -- & $3$ \\
\quad Actor hidden sizes & -- & $[155,127,105]$ \\
\quad Critic hidden sizes & -- & $[124,22,4]$ \\
\quad Initial log-std & $\log(\sigma_{\pi,0})$ & $0.0$ \\
\midrule
\multicolumn{3}{@{}l}{\textbf{PPO hyperparameters}} \\
\quad Initial learning rate & $\alpha_0$ & $2.0\times10^{-4}$ \\
\quad Final learning rate & $\alpha_{\mathrm{final}}$ & $1.0\times10^{-5}$ \\
\quad Discount factor & $\gamma$ & $0.9999$ \\
\quad GAE parameter & $\lambda_{\mathrm{GAE}}$ & $0.99$ \\
\quad Initial clip range & $\varepsilon_{\mathrm{clip},0}$ & $0.25$ \\
\quad Final clip range & $\varepsilon_{\mathrm{clip,final}}$ & $0.10$ \\
\quad Entropy coefficient & $c_{\mathrm{ent}}$ & $7.5\times10^{-4}$ \\
\quad Value function coefficient & $c_V$ & $0.6$ \\
\midrule
\multicolumn{3}{@{}l}{\textbf{Training configuration}} \\
\quad MC samples (training) & $N_s^{\mathrm{train}}$ & $512$ \\
\quad Parallel environments & $N_{\mathrm{env}}$ & $8$ \\
\quad Rollout length per environment & $n_{\mathrm{steps}}$ & $3200$ \\
\quad Number of minibatches & $N_{\mathrm{mb}}$ & $8$ \\
\quad PPO update epochs & $n_{\mathrm{epochs}}$ & $10$ \\
\quad Evaluation frequency & $n_{\mathrm{eval}}$ & $3200$ \\
\quad Total training timesteps & $N_{\mathrm{train}}$ & $8.0\times10^{7}$ \\
\bottomrule
\end{tabular}
\end{threeparttable}
\end{table}

\pagebreak

From a computational standpoint, the training performed offline using parallel environments required approximately 10 hours on a modern multi-core workstation.\footnote{Reported wall-clock time is implementation- and hardware-dependent.} At run time, policy evaluation requires only a single forward pass through the actor network at each decision node and is therefore computationally negligible relative to the underlying dynamics propagation and ensemble statistics.

\subsection{Results and Benchmark Evaluation}
\label{sec:em_results}

The Earth--Mars RCCRL policies are evaluated to address three key aspects of performance and robustness. First, the policy trained under Gaussian initial uncertainty is directly compared against the robust trajectory-optimization benchmark reported by Marmo et al.~\cite{marmo2025ChanceConstraintRobustTrajectory}. Second, the same RCCRL framework is evaluated under bounded uniform initial uncertainty with matched first and second moments to assess the effect of departing from the Gaussian regime. Third, the robustness of the learned controller to unmodeled process disturbances is quantified through a parametric noise-scaling study.

To make the Gaussian comparison meaningful, the results are reported against both the analytically evaluated and the corresponding MC validated values reported in~\cite{marmo2025ChanceConstraintRobustTrajectory}. In what follows, ``ROCP--CL'' denotes the reported analytical values, whereas ``ROCP--CL (MC)'' denotes the corresponding MC validation. Unless stated otherwise, all RCCRL statistics are computed from an evaluation campaign with $N_s^{\mathrm{eval}}=100{,}000$ independent realizations, matching the evaluation budget used in the ROCP--CL (MC) reference and reducing finite-sample variability in tail statistics.

\begin{table}[h!tb]
\centering
\caption{Performance comparison between RCCRL and ROCP--CL for the Earth--Mars transfer case study.}
\label{tab:em_results_comparison}
\begin{threeparttable}
\footnotesize
\begin{tabular}{@{}lccccc@{}}
\toprule
Quantity
& ROCP--CL~\cite{marmo2025ChanceConstraintRobustTrajectory}
& ROCP--CL (MC)~\cite{marmo2025ChanceConstraintRobustTrajectory}
& RCCRL (Gaussian)
& RCCRL (Uniform)
& Unit \\
\midrule
\multicolumn{6}{@{}l}{\textbf{Fuel cost}} \\
\quad $\Delta v_{\mathrm{tot},95}$                & $13.3725$                  & $13.3960$                  & $12.6594$                  & $13.0806$                  & km/s \\
\quad $\hat{\Delta v}_{\mathrm{tot}}$             & N/A                        & N/A                        & $11.9608$                  & $12.5911$                  & km/s \\
\quad $\Delta v_{\mathrm{nom}}$                   & $11.4918$                  & $11.4918$                  & $11.0861$                  & $11.7173$                  & km/s \\
\quad $\Delta v_{\mathrm{s},95}$\tnote{†}         & $1.8807$                   & $1.9042$                   & $1.5736$                   & $1.3633$                   & km/s \\
\quad $\hat{\Delta v}_{\mathrm{s}}$\tnote{††}      & N/A                        & N/A                        & $0.8746$                   & $0.8738$                   & km/s \\
\quad $\delta \Delta v_{\mathrm{tot},95}$\tnote{‡}
                                                  & N/A                        & $0.00$                     & $-0.7366$                  & $-0.3154$                  & km/s \\
\quad $\delta \Delta v_{\mathrm{tot},95}^{\%}$\tnote{‡‡}
                                                  & N/A                        & $0.00$                     & $-5.50$                    & $-2.35$                    & \% \\
\midrule
\multicolumn{6}{@{}l}{\textbf{Terminal dispersion}} \\
\quad $\sigma_{r_x,t_f}^{+}$                      & $1.3631\times10^{5}$       & $1.3130\times10^{5}$       & $7.6365\times10^{4}$       & $6.7153\times10^{4}$       & km \\
\quad $\sigma_{r_y,t_f}^{+}$                      & $1.4731\times10^{5}$       & $1.4663\times10^{5}$       & $3.2483\times10^{4}$       & $3.8833\times10^{4}$       & km \\
\quad $\sigma_{r_z,t_f}^{+}$                      & N/A                        & N/A                        & $6.8978\times10^{4}$       & $5.0647\times10^{4}$       & km \\
\quad $\sigma_{v_x,t_f}^{+}$                      & $8.5381\times10^{-3}$      & $7.5201\times10^{-3}$      & $5.2095\times10^{-3}$      & $6.9836\times10^{-3}$      & km/s \\
\quad $\sigma_{v_y,t_f}^{+}$                      & $9.3861\times10^{-3}$      & $8.5038\times10^{-3}$      & $5.3558\times10^{-3}$      & $5.0207\times10^{-3}$      & km/s \\
\quad $\sigma_{v_z,t_f}^{+}$                      & N/A                        & N/A                        & $4.4696\times10^{-3}$      & $5.6667\times10^{-3}$      & km/s \\
\midrule
\multicolumn{6}{@{}l}{\textbf{Terminal state at Mars}} \\
\quad $p_{\mathrm{SOI}}$                          & N/A                        & N/A                        & $1.00$                     & $1.00$                     & -- \\
\quad $\hat{q}_{0.95}(e_r)$                       & N/A                        & N/A                        & $2.2025\times10^{5}$       & $3.5492\times10^{5}$       & km \\
\quad $\min(e_r)$                                 & N/A                        & N/A                        & $7.0860\times10^{3}$       & $1.8378\times10^{5}$       & km \\
\quad $\max(e_r)$                                 & N/A                        & N/A                        & $3.9728\times10^{5}$       & $4.0450\times10^{5}$       & km \\
\quad $\E[e_r]$                                   & N/A                        & N/A                        & $1.4701\times10^{5}$       & $2.9785\times10^{5}$       & km \\
\quad $|\Delta \hat{r}_{x,t_f}^{+}|$              & N/A                        & $1.2135\times10^{2}$       & $1.3874\times10^{3}$       & $1.6248\times10^{3}$       & km \\
\quad $|\Delta \hat{r}_{y,t_f}^{+}|$              & N/A                        & $4.0769\times10^{2}$       & $1.0794\times10^{5}$       & $2.7287\times10^{5}$       & km \\
\quad $|\Delta \hat{r}_{z,t_f}^{+}|$              & N/A                        & N/A                        & $7.3792\times10^{3}$       & $8.3450\times10^{4}$       & km \\
\quad $|\Delta \hat{v}_{x,t_f}^{+}|$              & N/A                        & $6.9868\times10^{-6}$      & $1.4498\times10^{-4}$      & $6.1227\times10^{-5}$      & km/s \\
\quad $|\Delta \hat{v}_{y,t_f}^{+}|$              & N/A                        & $7.9084\times10^{-6}$      & $8.2600\times10^{-5}$      & $8.0702\times10^{-5}$      & km/s \\
\quad $|\Delta \hat{v}_{z,t_f}^{+}|$              & N/A                        & N/A                        & $1.2379\times10^{-4}$      & $5.7240\times10^{-5}$      & km/s \\
\bottomrule
\end{tabular}
\begin{tablenotes}
\footnotesize
\item[†] $\Delta v_{\mathrm{s},95} \doteq \Delta v_{\mathrm{tot},95}-\Delta v_{\mathrm{nom}}$.
\item[††] $\hat{\Delta v}_{\mathrm{s}} \doteq \hat{\Delta v}_{\mathrm{tot}}-\Delta v_{\mathrm{nom}}$.
\item[‡] $\delta \Delta v_{\mathrm{tot},95} = \Delta v_{\mathrm{tot},95}^{(j)} - \Delta v_{\mathrm{tot},95}^{\mathrm{ROCP\text{-}CL(MC)}}$.
\item[‡‡] $\delta \Delta v_{\mathrm{tot},95}^{\%}
=
100\,\delta \Delta v_{\mathrm{tot},95}
/\Delta v_{\mathrm{tot},95}^{\mathrm{ROCP\text{-}CL(MC)}}$.
\end{tablenotes}
\end{threeparttable}
\end{table}

\pagebreak

\tabref{tab:em_results_comparison} summarizes the key performance metrics for the two RCCRL policies and the ROCP--CL benchmarks. The primary fuel metric is the total impulse expenditure, summarized through both the empirical mean $\hat{\Delta v}_{\mathrm{tot}}$ and the 95th-percentile cost $\Delta v_{\mathrm{tot},95}$. Terminal performance is assessed through the terminal position error norm $e_r$ from \eqnref{eq:em_prob_capture}, the empirical Mars SOI-capture probability
\begin{equation}
p_{\mathrm{SOI}}
=
\frac{1}{N_s^{\mathrm{eval}}}
\sum_{i=1}^{N_s^{\mathrm{eval}}}
\mathbb{I}\!\left[e_r^{(i)} \le r_{\mathrm{SOI}}\right],
\label{eq:em_capture_prob}
\end{equation}
and the covariance-feasibility metric $\varepsilon_{\mathrm{cov}}$ defined in \eqnref{eq:em_cov_violation}. Together, these quantities characterize the cost of robustness, the probability of successful Mars capture, and the quality of the terminal state distribution.

\subsubsection{Gaussian Benchmark Case}
\label{subsec:em_results_gaussian}

The Gaussian benchmark case uses the initial-state distribution
$\mathcal{D}_0 = \mathcal{N}(\vect{x}_0,\mat{P}_0)$,
with the process-noise model of \secref{subsec:em_uncertainty} enabled during evaluation. This is the setting closest to the uncertainty assumptions used in the ROCP--CL benchmark and therefore provides the most direct comparison.

\begin{figure}[h!tb]
  \centering
  \includegraphics[width=0.98\textwidth]{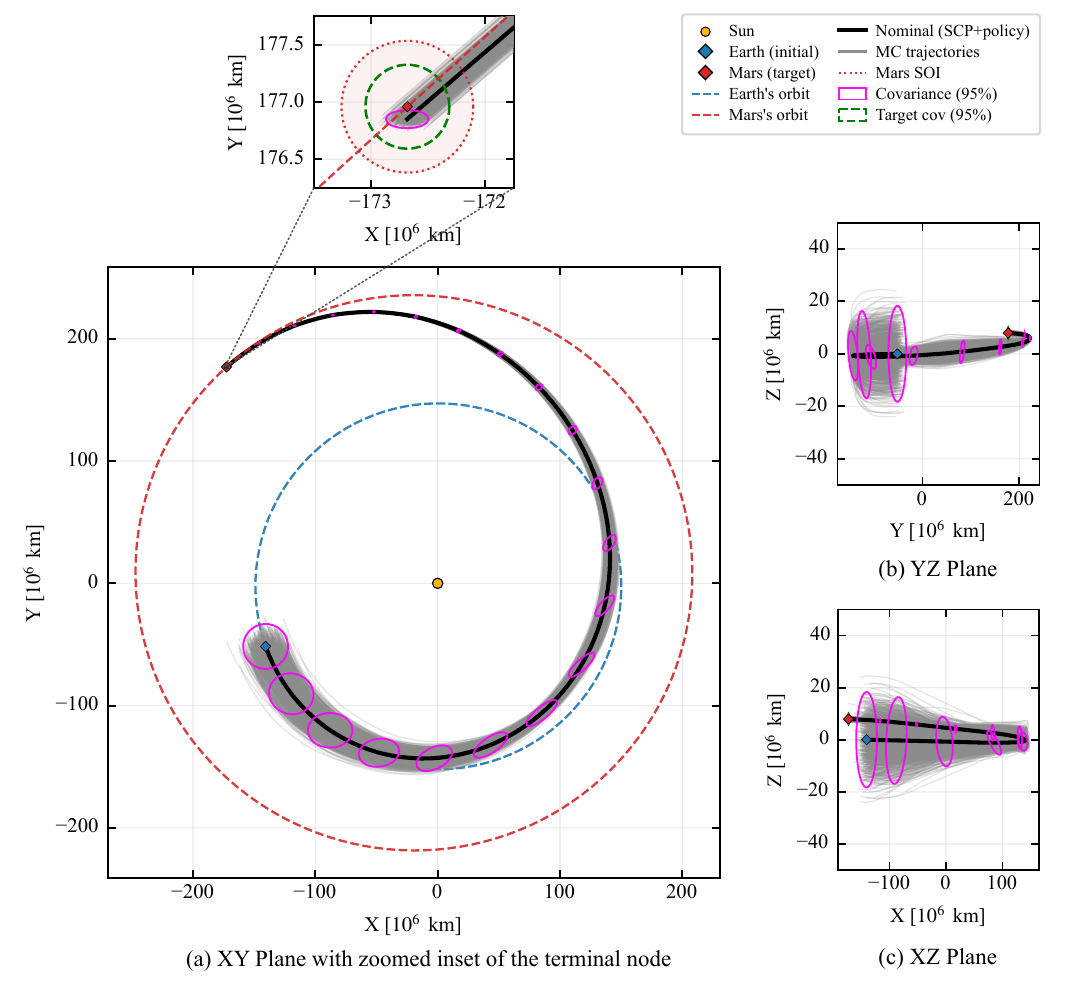}
  \caption{Closed-loop trajectory solution for the Earth--Mars transfer under Gaussian initial uncertainty. The trajectory bundle and covariance ellipses are magnified by a factor of $5$ for visualization, except in the terminal inset, which is shown at true scale. (a) In-plane heliocentric trajectory bundle with representative covariance ellipses and Mars-centered terminal inset. (b) $y$--$z$ projection. (c) $x$--$z$ projection.}
  \label{fig:em_gaussian_traj}
\end{figure}

\FloatBarrier

The Gaussian-case entries in \tabref{tab:em_results_comparison} show that the RCCRL policy attains a 95th-percentile total cost of
$\Delta v_{\mathrm{tot},95}=12.6594~\mathrm{km/s}$,
with empirical mean
$\hat{\Delta v}_{\mathrm{tot}}=11.9608~\mathrm{km/s}$.
Relative to the deterministic SCP baseline of \tabref{tab:em_scp_results}, the 95th-percentile cost increases by about $2.60~\mathrm{km/s}$ ($25.8\%$), reflecting the expected price of maintaining robustness under the imposed uncertainty and chance constraints. More importantly, relative to the ROCP--CL (MC) benchmark, RCCRL reduces the 95th-percentile total cost by about $0.737~\mathrm{km/s}$, corresponding to a $5.50\%$ reduction in the reported upper-tail total cost, while still maintaining empirical SOI capture and substantially tighter terminal dispersion. 

This reduction is accompanied by a lower reported nominal component $\Delta v_{\mathrm{nom}}$ and a lower proxy feedback contribution $\Delta v_{\mathrm{s},95}$ under the reporting convention of \tabref{tab:em_results_comparison}. 
\pagebreak
Although $\Delta v_{\mathrm{s}}$ is only a bookkeeping proxy and not a physically additive decomposition of the normed closed-loop control, the pair $\left(\Delta v_{\mathrm{nom}},\Delta v_{\mathrm{s},95}\right)$ is still informative as an indicator of how the policy reallocates control effort between systematic feedforward biasing and ensemble-dependent correction.

\begin{figure}[h!tb]
  \centering
  \includegraphics[width=0.98\textwidth]{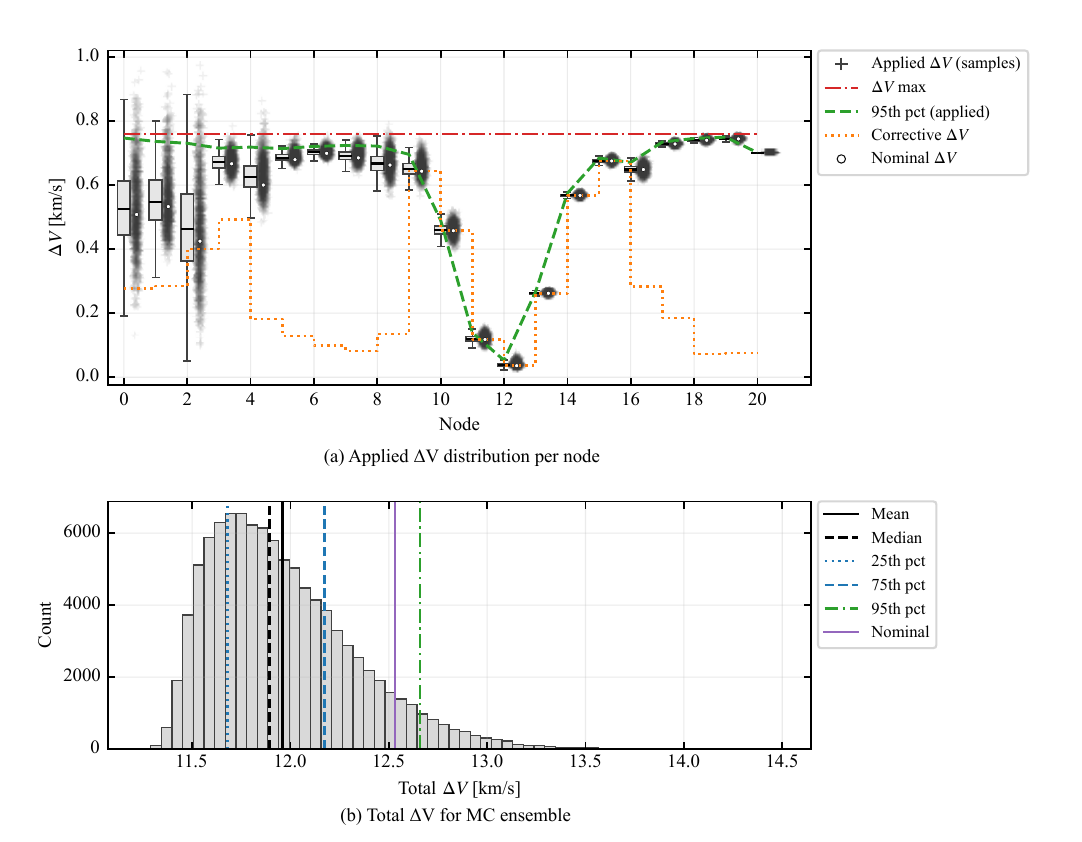}
  \caption{Control-effort statistics under Gaussian initial uncertainty. (a) Per-node impulse magnitude distribution across the evaluation ensemble. (b) Distribution of total closed-loop $\Delta v$ over the MC evaluation rollouts.}
  \label{fig:em_gaussian_ctrl}
\end{figure}

\FloatBarrier

\figref{fig:em_gaussian_traj} shows the resulting closed-loop trajectory ensemble. Across the transfer, the MC realizations remain tightly clustered around the nominal reference, with only modest out-of-plane dispersion. The terminal inset in \figref{fig:em_gaussian_traj}(a) shows that the final position cloud remains well inside the Mars SOI boundary, consistent with the empirical capture probability
$p_{\mathrm{SOI}} = 1.00$
over $100{,}000$ rollouts. This is an important result, showing that the reduced upper-tail fuel cost does not come at the expense of degraded capture reliability.

The control-effort statistics in \figref{fig:em_gaussian_ctrl} show that the largest spread in realized impulse magnitude appears early in the transfer. This is expected because the initial state dispersion is still prominent at the first maneuver nodes, and corrections applied early have the longest remaining horizon over which to influence downstream dispersion and terminal capture. Importantly, the 95th-percentile per-node impulse magnitude remains below the imposed bound $\Delta v_{\max}$ at every node, confirming satisfaction of the per-impulse chance constraint under the evaluation model with process noise.

\begin{figure}[h!tb]
  \centering
  \includegraphics[width=0.72\linewidth]{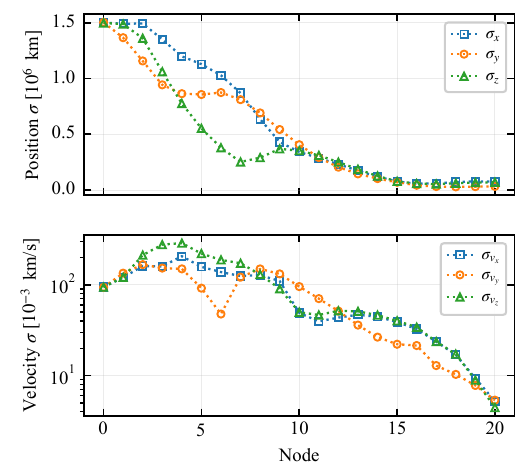}
  \caption{Evolution of ensemble dispersion along the trajectory for the Gaussian case. (a) Position-component standard deviations. (b) Velocity-component standard deviations on a logarithmic scale.}
  \label{fig:em_gaussian_sigma}
\end{figure}

\FloatBarrier

The evolution of ensemble dispersion is reported in \figref{fig:em_gaussian_sigma}. Position standard deviations contract by more than an order of magnitude over the transfer, reaching
$\sigma_{r_x,t_f}^{+}=7.64\times10^{4}~\mathrm{km}$,
$\sigma_{r_y,t_f}^{+}=3.25\times10^{4}~\mathrm{km}$, and
$\sigma_{r_z,t_f}^{+}=6.90\times10^{4}~\mathrm{km}$,
while the terminal velocity dispersions fall to the few-m/s level:
$\sigma_{v_x,t_f}^{+}=5.21\times10^{-3}~\mathrm{km/s}$,
$\sigma_{v_y,t_f}^{+}=5.36\times10^{-3}~\mathrm{km/s}$, and
$\sigma_{v_z,t_f}^{+}=4.47\times10^{-3}~\mathrm{km/s}$.
Compared with the ROCP--CL (MC) benchmark, the reported in-plane terminal dispersion is substantially tighter, with reductions of about $42\%$ in $\sigma_{r_x,t_f}^{+}$, $78\%$ in $\sigma_{r_y,t_f}^{+}$, $31\%$ in $\sigma_{v_x,t_f}^{+}$, and $37\%$ in $\sigma_{v_y,t_f}^{+}$.

\begin{figure}[h!tb]
  \centering
  \includegraphics[width=0.98\textwidth]{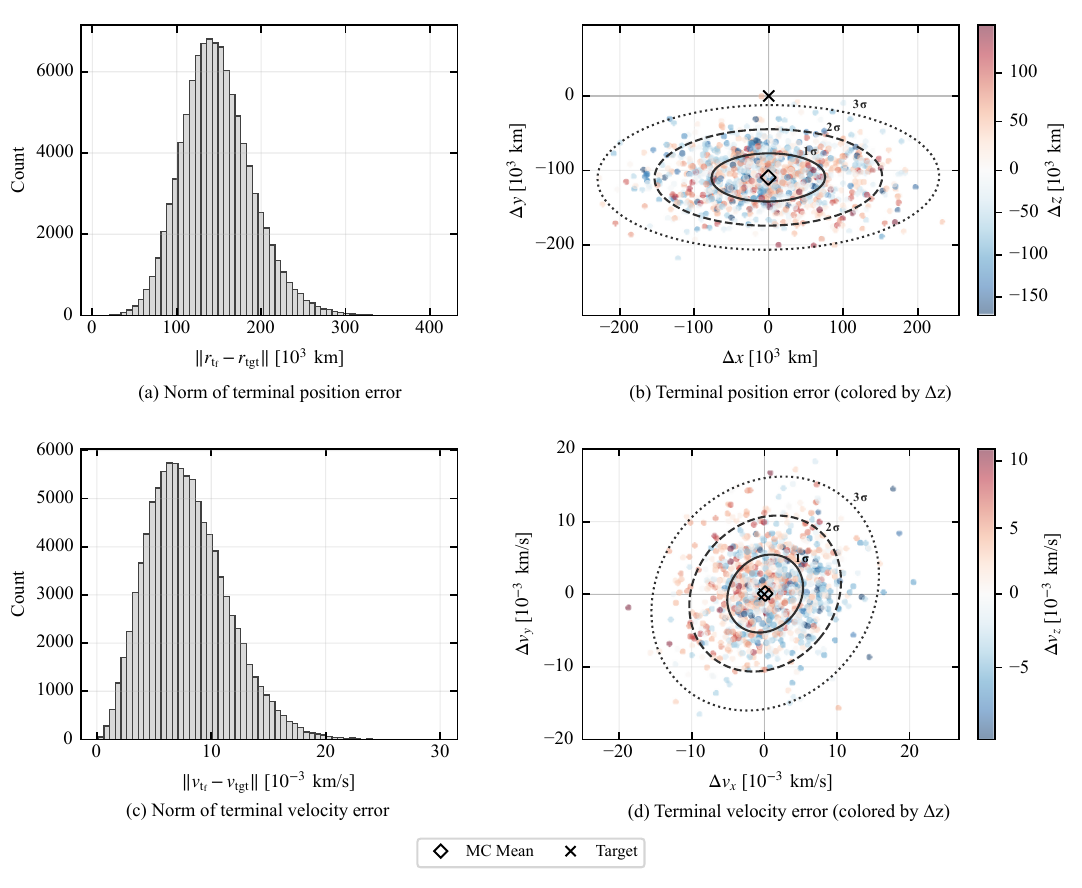}
  \caption{Terminal-state error distributions under Gaussian initial uncertainty. (a) Histogram of terminal position error norm $e_r$. (b) In-plane terminal position error scatter, colored by out-of-plane error. (c) Histogram of terminal velocity error norm. (d) In-plane terminal velocity error scatter, colored by out-of-plane error.}
  \label{fig:em_gaussian_disp}
\end{figure}

\figref{fig:em_gaussian_disp} characterizes the terminal-state distribution more directly. The 95th-percentile terminal position error is
$\hat{q}_{0.95}(e_r)=2.2025\times10^{5}~\mathrm{km}$,
and the worst observed realization is
$\max(e_r)=3.9728\times10^{5}~\mathrm{km}$,
both comfortably inside the Mars SOI radius of $5.77\times10^{5}~\mathrm{km}$. The terminal cloud also exhibits a non-zero mean bias,
$\E[e_r]=1.4701\times10^{5}~\mathrm{km}$,
which is expected because the objective prioritizes probabilistic capture, bounded upper-tail cost, and covariance feasibility, rather than explicit minimum-mean-miss objective. The covariance contraction alone does not prevent a coherent mean shift of the terminal distribution away from the target, which is precisely why the capture-oriented chance constraint is needed alongside the terminal covariance bound.

Taken together, the Gaussian-case results show that RCCRL can outperform the reported ROCP--CL benchmark in upper-tail total $\Delta v$ while preserving full capture feasibility and achieving tighter reported terminal dispersion. These results serve as the primary evidence that the learned affine correction law is competitive with an established chance-constrained robust trajectory-optimization method under comparable Gaussian uncertainty assumptions.

\subsubsection{Uniform Case}
\label{subsec:em_results_uniform}

To assess behavior beyond Gaussian assumptions, the uniform case evaluates the same RCCRL framework under the bounded initial-state model defined in \secref{subsec:em_uncertainty}. The distribution is constructed to match the Gaussian case in mean and marginal variance, but differs in support and higher-order statistics. This setting therefore tests whether the same chance-constraint machinery and affine policy structure remain effective when the uncertainty model changes without altering the first two moments.

\begin{figure}[h!tb]
  \centering
  \includegraphics[width=0.98\textwidth]{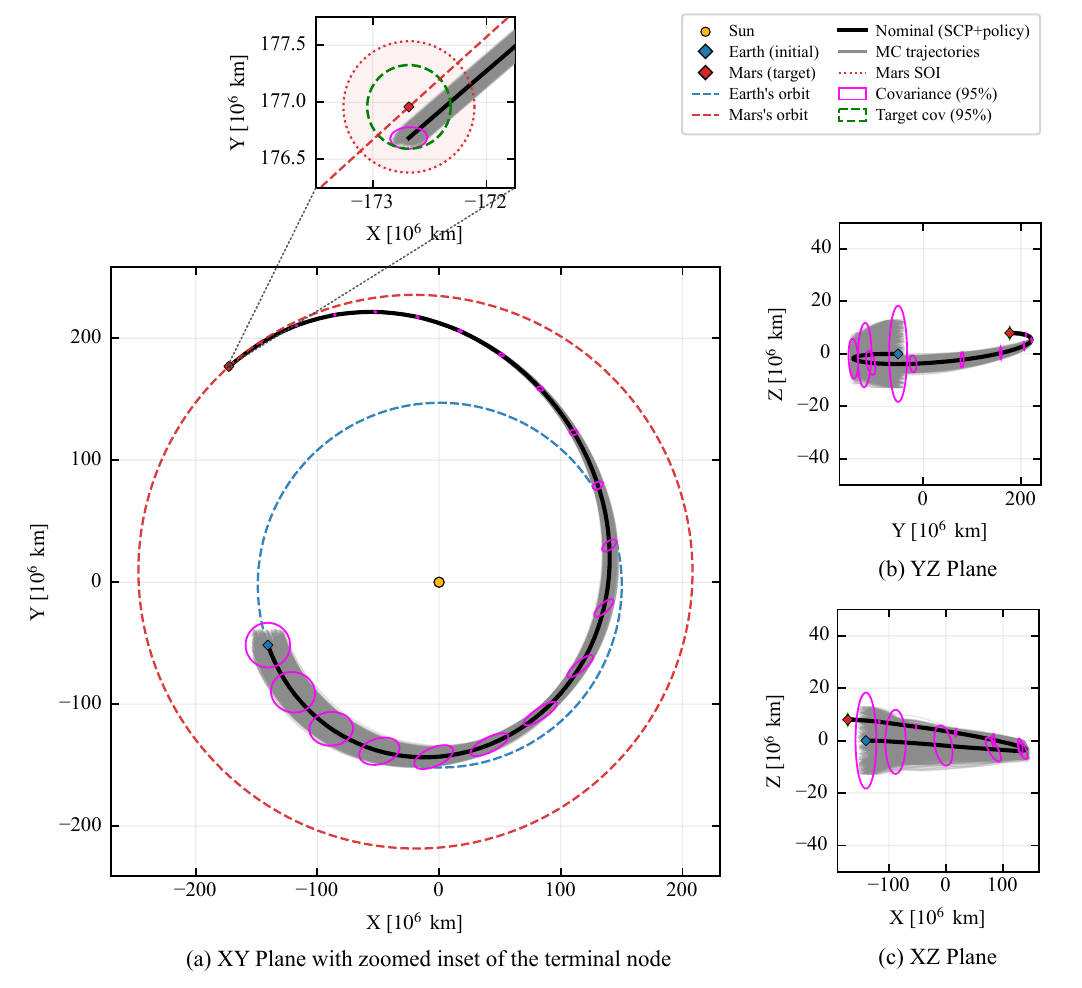}
  \caption{Closed-loop trajectory solution for the Earth--Mars transfer under uniform initial uncertainty. The trajectory bundle and covariance ellipses are magnified by a factor of $5$ for visualization, except in the terminal inset, which is shown at true scale. (a) In-plane heliocentric trajectory bundle with representative covariance ellipses and Mars-centered terminal inset. (b) $y$--$z$ projection. (c) $x$--$z$ projection.}
  \label{fig:em_uniform_traj}
\end{figure}

The uniform-case entries in the final column of \tabref{tab:em_results_comparison} show that, under uniform initial uncertainty, the learned policy attains a 95th-percentile total cost of
$\Delta v_{\mathrm{tot},95}=13.0806~\mathrm{km/s}$,
with empirical mean
$\hat{\Delta v}_{\mathrm{tot}}=12.5911~\mathrm{km/s}$.
Relative to the Gaussian RCCRL case, this corresponds to a modest increase of about $0.42~\mathrm{km/s}$ ($3.3\%$) in the 95th-percentile cost and $0.63~\mathrm{km/s}$ ($5.3\%$) in the mean cost. A plausible explanation is that, for fixed marginal variance, the bounded uniform model places more probability mass at relatively large initial deviations than a Gaussian distribution concentrated near the mean, so larger systematic corrective actions are needed to maintain feasibility under evaluation-time process noise.

\FloatBarrier

All probabilistic requirements remain satisfied in the uniform case, including the per-impulse bound, the SOI-capture condition, and terminal covariance feasibility. \figref{fig:em_uniform_traj} shows that the terminal cloud remains inside the SOI boundary, with
$p_{\mathrm{SOI}}=1.00$.
The corresponding control allocation shifts relative to the Gaussian case, such that the reported feedforward component increases to $\Delta v_{\mathrm{nom}}=11.7173~\mathrm{km/s}$, while the proxy upper-tail feedback contribution decreases to $\Delta v_{\mathrm{s},95}=1.3633~\mathrm{km/s}$. Although these quantities are only an interpretive decomposition, the redistribution is consistent with the policy relying more heavily on systematic feedforward biasing under bounded non-Gaussian uncertainty.

\begin{figure}[h!tb]
  \centering
  \includegraphics[width=0.98\textwidth]{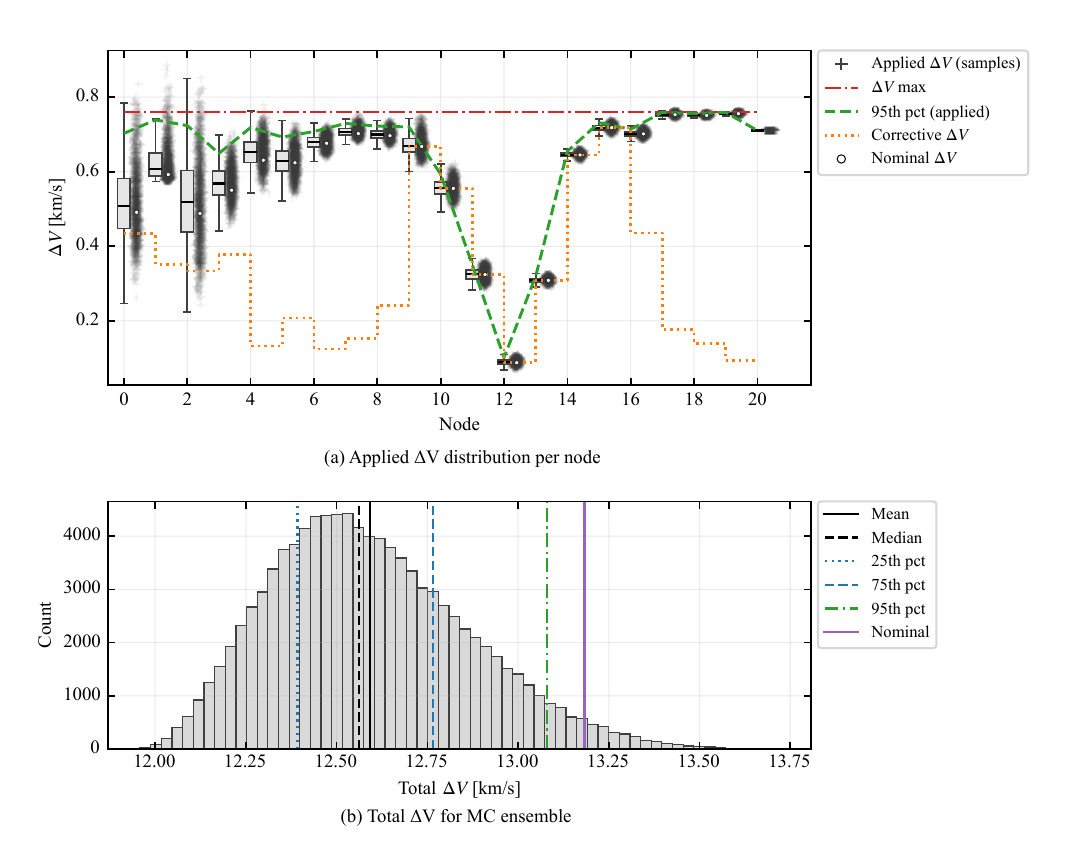}
  \caption{Control-effort statistics under uniform initial uncertainty. (a) Per-node impulse magnitude distribution across the evaluation ensemble. (b) Distribution of total closed-loop $\Delta v$ over the MC evaluation rollouts.}
  \label{fig:em_uniform_ctrl}
\end{figure}

\FloatBarrier

\begin{figure}[h!tb]
  \centering
  \includegraphics[width=0.72\linewidth]{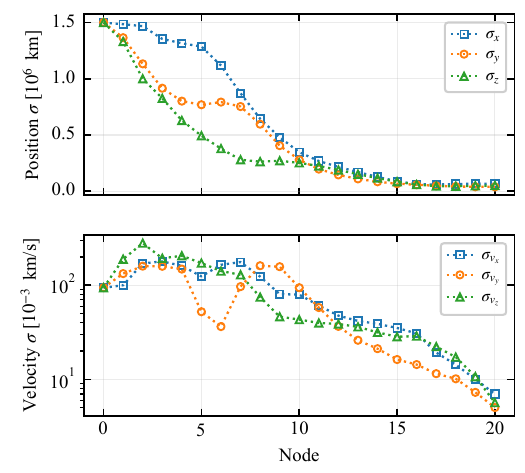}
  \caption{Evolution of ensemble dispersion along the trajectory for the uniform case. (a) Position-component standard deviations. (b) Velocity-component standard deviations on a logarithmic scale.}
  \label{fig:em_uniform_sigma}
\end{figure}

The control-effort and dispersion statistics in \figref{fig:em_uniform_ctrl} and \figref{fig:em_uniform_sigma} show the same qualitative structure as in the Gaussian case. The 95th-percentile per-node impulses remain below the bound, and the state dispersion contracts substantially toward the terminal epoch. At arrival, the position standard deviations are
$\sigma_{r_x,t_f}^{+}=6.72\times10^{4}~\mathrm{km}$,
$\sigma_{r_y,t_f}^{+}=3.88\times10^{4}~\mathrm{km}$, and
$\sigma_{r_z,t_f}^{+}=5.06\times10^{4}~\mathrm{km}$,
while the velocity standard deviations are
$\sigma_{v_x,t_f}^{+}=6.98\times10^{-3}~\mathrm{km/s}$,
$\sigma_{v_y,t_f}^{+}=5.02\times10^{-3}~\mathrm{km/s}$, and
$\sigma_{v_z,t_f}^{+}=5.67\times10^{-3}~\mathrm{km/s}$.
Compared with the Gaussian case, the terminal dispersion is redistributed across components rather than uniformly inflated, which is consistent with a different balance of corrective authority induced by the non-Gaussian initial-state distribution.

\begin{figure}[h!tb]
  \centering
  \includegraphics[width=0.98\textwidth]{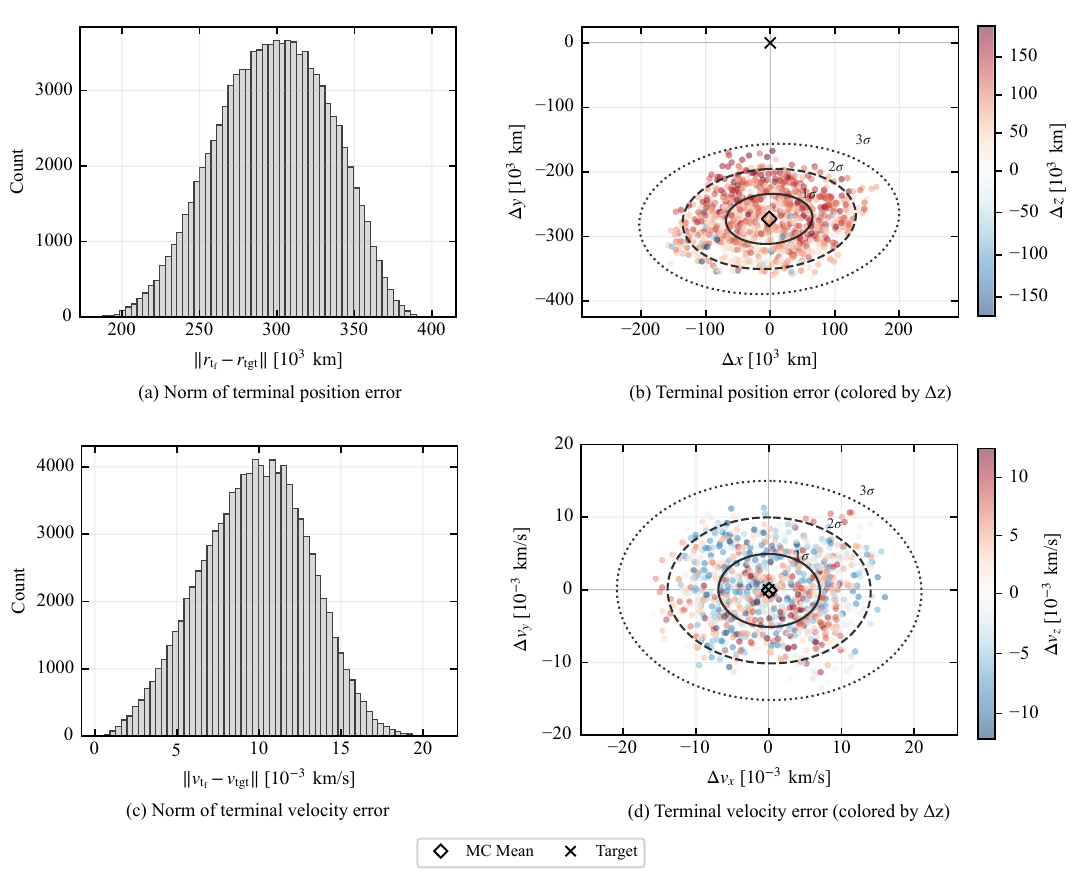}
  \caption{Terminal-state error distributions under uniform initial uncertainty. (a) Histogram of terminal position error norm $e_r$. (b) In-plane terminal position error scatter, colored by out-of-plane error. (c) Histogram of terminal velocity error norm. (d) In-plane terminal velocity error scatter, colored by out-of-plane error.}
  \label{fig:em_uniform_disp}
\end{figure}

The terminal position error distribution in \figref{fig:em_uniform_disp} yields
$\hat{q}_{0.95}(e_r)=3.5492\times10^{5}~\mathrm{km}$,
with
$\max(e_r)=4.0450\times10^{5}~\mathrm{km}$,
both still below the Mars SOI radius. The mean terminal miss distance increases to
$\E[e_r]=2.9785\times10^{5}~\mathrm{km}$,
which indicates a more pronounced coherent shift of the final distribution relative to the target state than in the Gaussian case. This again shows why matching only the first two moments of the initial uncertainty does not guarantee identical closed-loop terminal behavior under nonlinear propagation.

\pagebreak

Overall, the uniform-case results show that RCCRL remains feasible and cost-competitive under bounded non-Gaussian uncertainty without structural modification of the methodology. While the ROCP--CL benchmark is not designed for this setting and therefore cannot be compared on a like-for-like basis, the uniform-case 95th-percentile total cost still remains $0.3154~\mathrm{km/s}$ below the Gaussian ROCP--CL (MC) reference value, which provides a useful scale comparison.

\subsubsection{Robustness to Unmodeled Process Noise}
\label{subsec:em_results_noise}

As discussed in \secref{subsec:em_uncertainty}, process noise is disabled during training and reintroduced only at evaluation time. This simplifies training while creating a deliberate train--test mismatch. The purpose of the present subsection is to quantify how far the learned Gaussian-case policy can be pushed under unmodeled disturbances before feasibility degrades.

Let the nominal process-noise standard deviations be $(\sigma_{w,r},\sigma_{w,v})$. We introduce non-dimensional scaling factors $(\eta_r,\eta_v)$ and evaluate
\begin{equation}
(\sigma_{w,r},\sigma_{w,v})
\mapsto
(\eta_r \sigma_{w,r},\, \eta_v \sigma_{w,v}),
\qquad
\eta_r,\eta_v \ge 0,
\label{eq:em_noise_scaling_results}
\end{equation}
such that $(\eta_r,\eta_v)=(1,1)$ reproduces the nominal evaluation noise level used in the Gaussian and uniform results. Because the disturbance covariance scales quadratically in the standard deviations, the discrete process-noise covariance obeys
\begin{equation}
\mat{Q}(\eta_r,\eta_v)
\propto
\diag(\eta_r^2,\eta_v^2).
\label{eq:em_noise_cov_scaling}
\end{equation}

\FloatBarrier

For each $(\eta_r,\eta_v)$, the following quantities are monitored: (i) the 95th-percentile total cost $\Delta v_{\mathrm{tot},0.95}$, (ii) the Mars capture probability $p_{\mathrm{SOI}}$, (iii) the terminal-distance quantile $\hat{q}_{0.95}(e_r)$, (iv) the covariance-feasibility violation $\varepsilon_{\mathrm{cov}}$, and (v) the feasibility indicator
\begin{equation}
\mathcal{F}(\eta_r,\eta_v)
=
\mathbb{I}\Big[
\max_k \hat{q}_{0.95}(\norm{\Delta \vect{v}_k}) \le \Delta v_{\max}
\;\wedge\;
\hat{q}_{0.95}(e_r) \le r_{\mathrm{SOI}}
\;\wedge\;
\varepsilon_{\mathrm{cov}} = 0
\Big].
\label{eq:em_feasibility_indicator}
\end{equation}
This is interpreted in a feasibility-first sense, i.e., cost trends are only meaningful while $\mathcal{F}=1$.

\begin{figure}[h!tb]
  \centering
  \includegraphics[width=0.98\textwidth]{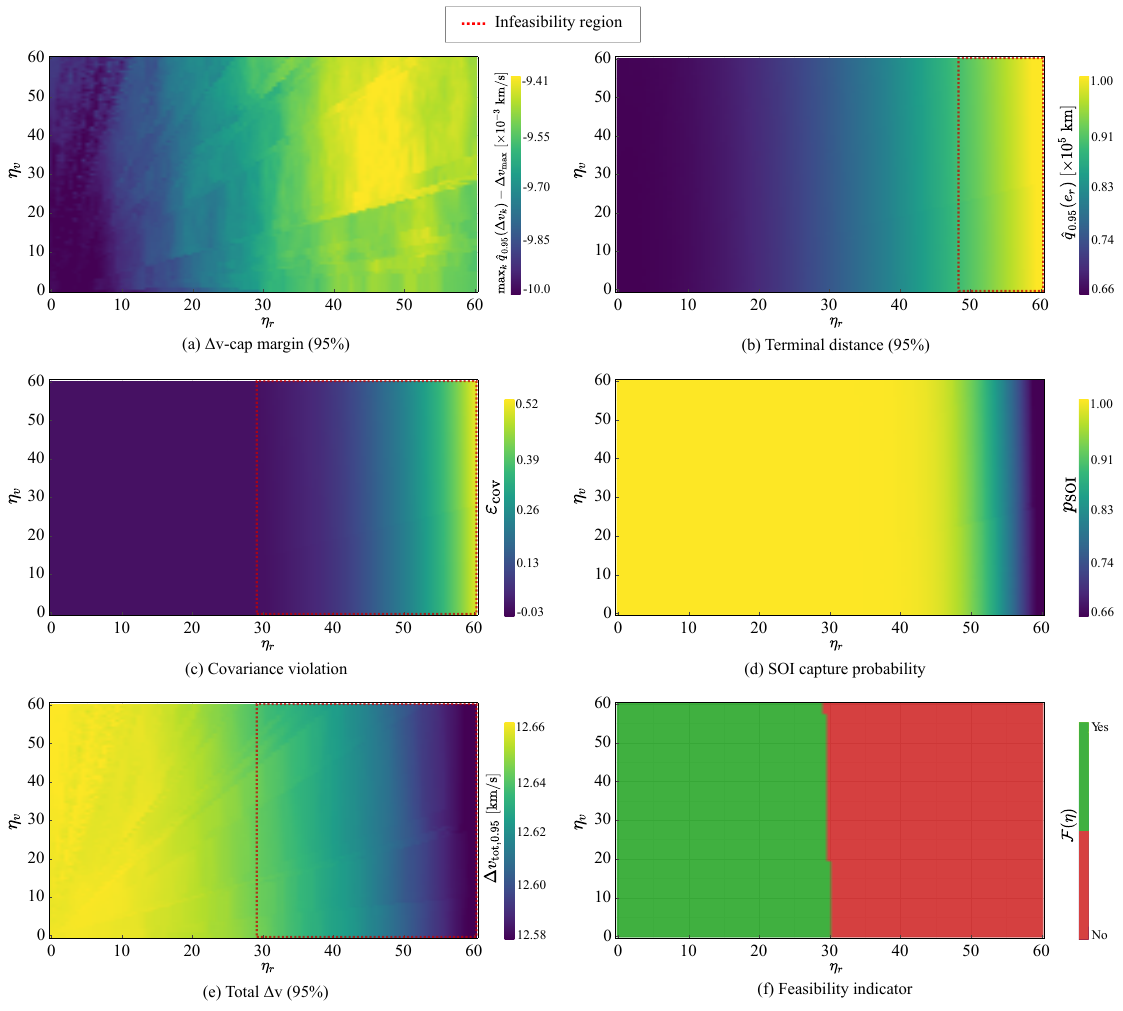}
  \caption{Robustness to two-dimensional variation in unmodeled process noise under Gaussian initial uncertainty. Independent scalings $(\eta_r,\eta_v)$ are applied to the position and velocity disturbance standard deviations. Panels show: (a) worst-node $\Delta v$-cap margin, (b) terminal-distance quantile $\hat q_{0.95}(e_r)$, (c) covariance-constraint violation $\varepsilon_{\mathrm{cov}}$, (d) Mars SOI capture probability $p_{\mathrm{SOI}}$, (e) total cost $\Delta v_{\mathrm{tot},0.95}$, and (f) feasibility indicator $\mathcal{F}(\eta_r,\eta_v)$.}
  \label{fig:em_noise_heatmap}
\end{figure}

\figref{fig:em_noise_heatmap} reports the two-dimensional sweep over $(\eta_r,\eta_v)$ and reveals several aspects of how the trained RCCRL policy responds to disturbances not seen during training. The per-impulse $\Delta v$ cap remains inactive over the explored range, indicating that loss of robustness is not driven by control-authority saturation. Instead, feasibility is first lost through the covariance constraint, as $\varepsilon_{\mathrm{cov}}$ becomes non-zero before the SOI-capture condition fails. The robustness boundary is visibly more sensitive to the position-noise scaling $\eta_r$ than to the velocity-noise scaling $\eta_v$, indicating that positional diffusion between decision nodes is the dominant driver of constraint loss for this transfer and update schedule.

\tabref{tab:em_noise_sensitivity} provides a compact one-dimensional summary along the diagonal slice $\eta_r=\eta_v=\eta$, which preserves the nominal ratio between position and velocity disturbance intensities. Up to about $\eta=15$, the policy performance remains essentially unchanged: both $\Delta v_{\mathrm{tot},0.95}$ and $\hat{\Delta v}_{\mathrm{tot}}$ vary only weakly, $p_{\mathrm{SOI}}$ remains unity, and $\varepsilon_{\mathrm{cov}}=0$. The first infeasible point appears at
$\eta_{\mathrm{crit}} \approx 29.7$,
where feasibility is lost through covariance violation while the terminal-distance quantile remains comfortably within the SOI radius and the capture probability is still $1.00$.

At higher noise levels, the terminal-distance quantile continues to grow and approaches the SOI radius only later, around $\eta \approx 48.9$, at which point the capture probability begins to degrade. By $\eta=60$, the capture probability drops to $61.4\%$. This ordering strongly indicates that the dominant robustness limitation in this case study is not immediate control saturation, but rather the growth of terminal dispersion beyond the prescribed covariance-feasibility envelope.

\begin{table}[h!tb]
\centering
\caption{Noise-sensitivity summary for the Gaussian RCCRL policy under scaled process noise.}
\label{tab:em_noise_sensitivity}
\begin{threeparttable}
\footnotesize
\begin{tabular}{@{}lcccccc@{}}
\toprule
$\eta$ & $0.0$ & $1.0$ & $15.0$ & $29.7$ & $48.9$ & $60.0$ \\
Regime & Baseline & Nominal & Moderate & $\eta_{\mathrm{crit}}$ & SOI-limit & Max tested \\
\midrule
$\Delta v_{\mathrm{tot},0.95}$ [km/s] & $12.6595$ & $12.6594$ & $12.6534$ & $12.6360$ & $12.5984$ & $12.5681$ \\
$\hat{\Delta v}_{\mathrm{tot}}$ [km/s] & $11.9601$ & $11.9608$ & $11.9543$ & $11.9375$ & $11.8992$ & $11.8689$ \\
$p_{\mathrm{SOI}}$ [--]                & $1.000$   & $1.000$   & $1.000$   & $1.000$   & $0.949$   & $0.614$   \\
$\hat{q}_{0.95}(e_r)$ [$\times 10^5$ km]
                                        & $2.2124$  & $2.2025$  & $2.5539$  & $3.5735$  & $5.7763$  & $7.4237$  \\
$\varepsilon_{\mathrm{cov}}$ [--]       & $0$       & $0$       & $0$       & $3.7929\times10^{-5}$ & $1.9869\times10^{-1}$ & $5.1051\times10^{-1}$ \\
Feasible                                & Yes       & Yes       & Yes       & No        & No        & No \\
\bottomrule
\end{tabular}
\end{threeparttable}
\end{table}

\begin{figure}[h!tb]
  \centering
  \includegraphics[width=0.72\linewidth]{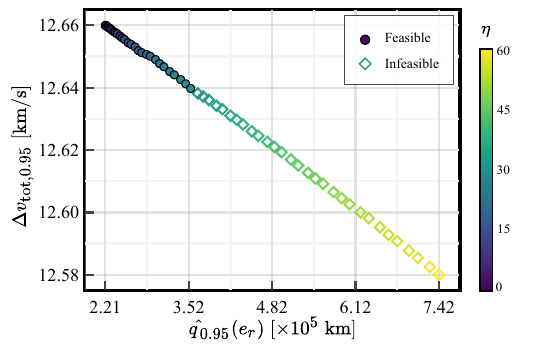}
  \caption{Trade-off between terminal dispersion $\hat q_{0.95}(e_r)$ and 95th-percentile total cost $\Delta v_{\mathrm{tot},0.95}$ over the diagonal sweep $\eta_r=\eta_v=\eta$. Marker style indicates feasibility and color encodes the process-noise scaling level.}
  \label{fig:em_noise_pareto}
\end{figure}

The control-cost--terminal-dispersion trade-off is summarized in \figref{fig:em_noise_pareto}. Although $\Delta v_{\mathrm{tot},0.95}$ decreases mildly as $\eta$ increases, this trend occurs predominantly after feasibility has already been lost and should therefore not be misread. Rather, once the covariance and eventually capture constraints are no longer maintained, the controller expends less corrective effort while allowing terminal dispersion to grow.

The process-noise study shows that the learned closed-loop structure provides meaningful robustness to moderate unmodeled disturbances without retraining, but with a finite robustness radius determined by disturbance-driven dispersion growth between decision nodes relative to the corrective authority learned during training. If stronger guarantees are required under substantially larger disturbance levels, that robustness radius can be further extended either by including process noise during training or by increasing allowable feedback authority, with an associated fuel-cost trade-off.

Taken together, the Earth--Mars results establish the primary benchmark case of the paper. Under Gaussian uncertainty, RCCRL is competitive with an established chance-constrained robust trajectory-optimization reference in upper-tail fuel cost while preserving full capture feasibility. 
\pagebreak
Under bounded uniform uncertainty and under evaluation-time process-noise mismatch, the same framework remains feasible and interpretable without structural modification, indicating that its usefulness is not confined to a single uncertainty family or an exactly matched training environment.
\section{Case Study~\Romannum{2}: Atmospheric Pinpoint Rocket Landing}
\label{sec:cs2_rktland}

The second case study is used to test the portability of RCCRL beyond the long-horizon impulsive-transfer regime. Unlike the months-long impulsive interplanetary transfer of the Earth--Mars benchmark in Case Study~\Romannum{1}, the present problem is a short-horizon continuous-thrust powered descent, which unfolds over $20$~s, with mass depletion, aerodynamic drag, and state-dependent glide-slope and thrust-magnitude constraints. The same RCCRL framework is instantiated with only problem-specific changes in dynamics, scaling, and reward metrics, and the intended purpose is to assess whether the same nominal-backbone, affine-correction, and sample-based chance-constraint scaffold remains effective when carried into a materially different trajectory design problem.

\subsection{Problem Setup and RCCRL Instantiation}
\label{subsec:rktland_setup}

We consider a two-dimensional atmospheric pinpoint landing problem in a local tangent-plane frame centered on the touchdown site. The translational state is
\begin{equation}
\vect{x}(t)=
\begin{bmatrix}
x(t)\\
y(t)\\
v_x(t)\\
v_y(t)
\end{bmatrix}
\in \R^{4},
\label{eq:rl_state}
\end{equation}
where $x$ is horizontal position, $y$ is altitude, and $(v_x,v_y)$ are the corresponding velocity components. The vehicle mass $m(t)$ is propagated alongside the state.

The continuous-time dynamics are
\begin{subequations}
\label{eq:rl_dynamics}
\begin{align}
\dot{x}   &= v_x, \qquad \dot{y}=v_y,
\\
\dot{v}_x &= -\frac{\Gamma_d}{m}\norm{\vect{v}}\,v_x + U_x,
\\
\dot{v}_y &= -g - \frac{\Gamma_d}{m}\norm{\vect{v}}\,v_y + U_y,
\\
\dot{m}   &= -\frac{m\norm{\vect{U}}}{I_{\mathrm{sp}}g_0},
\end{align}
\end{subequations}
where $g$ is constant gravitational acceleration, $\Gamma_d=\tfrac{1}{2}\rho_0 C_d S_{\mathrm{ref}}$ is the drag parameter combining sea-level density $\rho_0$, drag coefficient $C_d$, and reference area $S_{\mathrm{ref}}$, and $\vect{U}=[U_x,U_y]^\top$ is the acceleration due to thrust. The vehicle begins from the initial state and mass listed in \tabref{tab:rl_problem_data} and must achieve a soft landing at
\begin{equation}
\vect{x}(t_f)=\vect{0}.
\label{eq:rl_terminal_target}
\end{equation}

The descent is discretized over $N$ uniform segments with time grid
\begin{equation}
t_k=t_0+k\Delta t,
\qquad
k=0,\dots,N,
\qquad
\Delta t=\frac{\Delta T}{N}.
\label{eq:rl_time_grid}
\end{equation}
Unlike the Earth--Mars transcription, the control is held constant over each segment and the state is advanced by numerical integration of \eqnref{eq:rl_dynamics}. The deterministic objective is to maximize final mass, equivalently minimizing propellant use.

\pagebreak

Two path constraints are imposed. First, the thrust magnitude is limited by the available engine authority:
\begin{equation}
\norm{\vect{U}_k}\le \frac{T_{\max}}{m_k},
\qquad
k=0,\dots,N-1.
\label{eq:rl_thrust_constraint}
\end{equation}
Second, the descent must remain within a conic glide-slope corridor:
\begin{equation}
|x_k|\le y_k\tan(\theta_{\max}),
\qquad
k=0,\dots,N-1.
\label{eq:rl_glideslope_constraint}
\end{equation}

As in the Earth--Mars case, all variables are non-dimensionalized to improve numerical conditioning. The adopted physical parameters, boundary conditions, uncertainty levels, and discretization settings are summarized in \tabref{tab:rl_problem_data}.

\begin{table}[h!tb]
\centering
\caption{Problem parameters for the atmospheric rocket landing case study.}
\label{tab:rl_problem_data}
\begin{threeparttable}
\footnotesize
\begin{tabular}{@{}lcc@{}}
\toprule
Parameter & Value & Unit \\
\midrule
\multicolumn{3}{@{}l}{\textbf{Normalization reference}} \\
\quad Reference length scale, $r_{\mathrm{scale}}$ & $3{,}000$ & m \\
\quad Reference acceleration, $a_{\mathrm{scale}}$ & $9.81$ & m/s$^2$ \\
\quad Reference mass, $m_{\mathrm{scale}}$ & $55{,}000$ & kg \\
\quad Reference velocity scale, $v_{\mathrm{scale}}$ & $171.55$ & m/s \\
\quad Reference time scale, $t_{\mathrm{scale}}$ & $17.49$ & s \\
\midrule
\multicolumn{3}{@{}l}{\textbf{Discretization and problem setup}} \\
\quad Total time of flight, $\Delta T$ & $20.0$ & s \\
\quad Number of segments, $N$ & $40$ & -- \\
\quad Segment duration, $\Delta t$ & $0.5$ & s \\
\quad Maximum glide-slope angle, $\theta_{\max}$ & $70$ & deg \\
\quad Risk level, $\beta$ & $5\%$ & -- \\
\quad Process-noise intensity, $g_v$ & $0.25$ & m/s$^{3/2}$ \\
\midrule
\multicolumn{3}{@{}l}{\textbf{Vehicle and environment parameters}} \\
\quad Maximum thrust, $T_{\max}$ & $1{,}375.60$ & kN \\
\quad Specific impulse, $I_{\mathrm{sp}}$ & $443$ & s \\
\quad Standard gravity, $g_0$ & $9.81$ & m/s$^2$ \\
\quad Sea-level density, $\rho_0$ & $1.225$ & kg/m$^3$ \\
\quad Drag coefficient, $C_d$ & $0.5$ & -- \\
\quad Reference area, $S_{\mathrm{ref}}$ & $12.54$ & m$^2$ \\
\midrule
\multicolumn{3}{@{}l}{\textbf{Initial and target states}} \\
\quad Initial position, $[x_{t_0},y_{t_0}]^\top$ & $[950,\,3{,}000]^\top$ & m \\
\quad Initial velocity, $[v_{x,t_0},v_{y,t_0}]^\top$ & $[-118.33,\,-231.51]^\top$ & m/s \\
\quad Initial mass, $m_{t_0}$ & $55{,}000$ & kg \\
\quad Target position, $[x_{t_f},y_{t_f}]^\top$ & $[0,\,0]^\top$ & m \\
\quad Target velocity, $[v_{x,t_f},v_{y,t_f}]^\top$ & $[0,\,0]^\top$ & m/s \\
\midrule
\multicolumn{3}{@{}l}{\textbf{Initial and target uncertainty}} \\
\quad Initial position std. dev., $\sigma_{r,t_0}$ & $10.00$ & m \\
\quad Initial velocity std. dev., $\sigma_{v,t_0}$ & $3.1623$ & m/s \\
\quad Target position std. dev., $\sigma_{r,t_f}$ & $1.00$ & m \\
\quad Target velocity std. dev., $\sigma_{v,t_f}$ & $1.00$ & m/s \\
\bottomrule
\end{tabular}
\end{threeparttable}
\end{table}

The uncertain initial state is modeled either as a Gaussian distribution (\eqnref{eq:em_initial_gaussian}) or as a bounded uniform distribution with matched first and second moments, as in the Earth--Mars case. A continuous-time process-noise model is also included through additive acceleration-level diffusion and discretized segment-by-segment during evaluation. In the nominal experiments reported here, process noise is primarily used at evaluation time in order to probe out-of-distribution robustness rather than being treated as an always-on training disturbance.

\pagebreak

The landing problem includes three categories of constraints. First, the thrust magnitude must satisfy a node-wise probabilistic bound corresponding to \eqnref{eq:rl_thrust_constraint}. Second, the glide-slope corridor is enforced on the ensemble mean trajectory,
\begin{equation}
|\hat{x}_k| \le \hat{y}_k \tan(\theta_{\max}),
\qquad
k=0,\dots,N-1,
\label{eq:rl_mean_glideslope}
\end{equation}
and is therefore used here as a nominal descent-geometry constraint rather than as a sample-wise terrain-clearance probabilistic constraint. This relaxation avoids over-constraining individual dispersed samples while the terminal distribution is controlled through covariance-based penalties. For stricter mission requirements, a quantile penalty could instead be imposed on the sample-wise residual without changing the underlying RCCRL framework. Third, terminal dispersion is monitored through the same eigenvalue-based covariance-feasibility violation metric introduced in \secref{subsec:em_uncertainty}.

The deterministic nominal backbone is generated offline with SCP, using a thrust-continuation strategy to improve convergence. The resulting fuel-optimal baseline achieves sub-centimeter terminal accuracy with final mass $50{,}105.85~\mathrm{kg}$. This corresponds to propellant consumption of $4{,}894.15~\mathrm{kg}$ and equivalent velocity increment $\Delta v_{\mathrm{eq}}=405.01~\mathrm{m/s}$. This nominal solution provides the reference state and control histories around which RCCRL learns its robust closed-loop augmentation.

For this landing problem, the RCCRL action comprises a feedforward thrust correction $\vect{U}^{\mathrm{corr}}_k\in\R^2$ and a time-varying gain matrix $\mat{K}_k\in\R^{2\times4}$:
\begin{equation}
\vect{a}_k=
\begin{bmatrix}
\vect{U}^{\mathrm{corr}}_k \\
\mathrm{vec}(\mat{K}_k)
\end{bmatrix}.
\label{eq:rl_action}
\end{equation}
The applied control for the $i$-th ensemble member is
\begin{equation}
\vect{U}^{(i)}_k
=
\vect{U}^{\mathrm{ref}}_k
+
\vect{U}^{\mathrm{corr}}_k
+
\mat{K}_k\!\left(\vect{x}^{(i)}_k-\hat{\vect{x}}_k\right),
\label{eq:rl_ctrl_law}
\end{equation}
where $\vect{U}^{\mathrm{ref}}_k$ is the nominal SCP control and $\hat{\vect{x}}_k$ is the empirical ensemble mean. The action space therefore has dimension $10$.

\begin{table}[h!tb]
\centering
\caption{RCCRL reward parameters for the rocket landing case study.}
\label{tab:rl_reward_config}
\begin{threeparttable}
\footnotesize
\begin{tabular}{@{}lcc@{}}
\toprule
Parameter & Symbol & Value \\
\midrule
\multicolumn{3}{@{}l}{\textbf{Reward weights}} \\
\quad Mass consumption weight & $\lambda_m$ & $0.1$ \\
\quad Thrust violation weight & $\lambda_T$ & $0.01$ \\
\quad Glide-slope violation weight & $\lambda_\theta$ & $1.0$ \\
\quad Terminal state weight & $\lambda_{\mathrm{state}}$ & $20.0$ \\
\quad Position error weight & $\lambda_r$ & $2.0$ \\
\quad Velocity error weight & $\lambda_v$ & $15.0$ \\
\quad Covariance violation weight & $\lambda_{\mathrm{cov}}$ & $100.0$ \\
\quad Terminal bonus & -- & $200.0$ \\
\midrule
\multicolumn{3}{@{}l}{\textbf{Terminal tolerances}} \\
\quad State error threshold & $\varepsilon_{\mathrm{state}}^{\max}$ & $15.0$ \\
\quad Covariance violation threshold & $\varepsilon_{\mathrm{cov}}^{\max}$ & $0.5$ \\
\bottomrule
\end{tabular}
\end{threeparttable}
\end{table}

The observation includes the empirical mean state, mean mass, half-vectorized covariance, reference control, and normalized time-to-go, yielding an $18$-dimensional input. As in the Earth--Mars case, the policy is trained with PPO using the generic actor--critic structure of \secref{subsec:framework_policy}, resized to the reduced input and output dimensions, leading to the corresponding hidden layer sizes of $[90,67,50]$ and $[72,16,4]$ for the actor and critic networks, respectively. The training ensemble size is $512$, and the total number of training timesteps is $8.0\times10^7$. The same hyperparameters are used for both the Gaussian and uniform cases, with no tuning or normalization changes between them. From a computational standpoint, training for each landing case is performed offline using 8 parallel environments and required approximately 6 hours on a modern multi-core workstation.\footnote{Reported wall-clock time is implementation- and hardware-dependent.}

The reward follows the same generic RCCRL structure as in \secref{subsec:framework_objective}, but specializes the performance and feasibility metrics to the powered-descent setting. Frequent (per-step) penalties account for propellant consumption, thrust-cap violation, and mean glide-slope violation, while terminal penalties act on touchdown state error and covariance infeasibility. The adopted weights are listed in \tabref{tab:rl_reward_config}.

\subsection{Results Under Gaussian and Uniform Uncertainty}
\label{subsec:rktland_results}

Unless stated otherwise, all reported statistics are obtained from a Monte Carlo evaluation campaign with $N_s^{\mathrm{eval}}=100{,}000$ independent realizations. \tabref{tab:rl_results_comparison} summarizes the key landing metrics for the Gaussian and uniform policies, both without and with nominal process noise enabled during evaluation.

\begin{table*}[h!tb]
\centering
\caption{Performance comparison of RCCRL policies for the rocket landing case study.}
\label{tab:rl_results_comparison}
\begin{threeparttable}
\footnotesize
\begin{tabular}{@{}lccccc@{}}
\toprule
Quantity & RCCRL (Gaussian) & RCCRL (Gaussian) & RCCRL (Uniform) & RCCRL (Uniform) & Unit \\
         & w/o noise        & w/ noise         & w/o noise       & w/ noise        &      \\
\midrule
\multicolumn{6}{@{}l}{\textbf{Final mass}} \\
\quad $\hat{m}_{f}$               & $50{,}077.77$ & $50{,}077.95$ & $50{,}071.76$ & $50{,}071.66$ & kg \\
\quad $m_{f,05}$                  & $50{,}015.71$ & $50{,}015.51$ & $50{,}011.35$ & $50{,}008.62$ & kg \\
\quad $m_{f,95}$                  & $50{,}131.74$ & $50{,}132.98$ & $50{,}134.40$ & $50{,}136.07$ & kg \\
\midrule
\multicolumn{6}{@{}l}{\textbf{Propellant consumption}} \\
\quad $\Delta \hat{m}$            & $4{,}922.23$  & $4{,}922.05$  & $4{,}928.24$  & $4{,}928.34$  & kg \\
\quad $\Delta m_{95}$             & $4{,}984.29$  & $4{,}984.49$  & $4{,}988.65$  & $4{,}991.38$  & kg \\
\quad $\Delta m_{05}$             & $4{,}868.26$  & $4{,}867.02$  & $4{,}865.60$  & $4{,}863.93$  & kg \\
\midrule
\multicolumn{6}{@{}l}{\textbf{Equivalent $\Delta v$}} \\
\quad $\Delta \hat{v}_{\mathrm{eq}}$ & $407.45$ & $407.43$ & $407.97$ & $407.98$ & m/s \\
\quad $\Delta v_{\mathrm{eq},95}$    & $412.84$ & $412.85$ & $413.22$ & $413.45$ & m/s \\
\quad $\Delta v_{\mathrm{eq},05}$    & $402.77$ & $402.66$ & $402.54$ & $402.39$ & m/s \\
\midrule
\multicolumn{6}{@{}l}{\textbf{Terminal state dispersion}} \\
\quad $|\Delta \hat{r}_{x,t_f}|$  & $1.76\times10^{-2}$ & $4.89\times10^{-2}$ & $1.95\times10^{-2}$ & $4.27\times10^{-2}$ & m \\
\quad $|\Delta \hat{r}_{y,t_f}|$  & $5.22\times10^{-2}$ & $1.81\times10^{-1}$ & $1.10\times10^{-1}$ & $1.89\times10^{-1}$ & m \\
\quad $|\Delta \hat{v}_{x,t_f}|$  & $4.07\times10^{-3}$ & $1.02\times10^{-3}$ & $7.73\times10^{-3}$ & $4.43\times10^{-3}$ & m/s \\
\quad $|\Delta \hat{v}_{y,t_f}|$  & $2.54\times10^{-3}$ & $6.09\times10^{-3}$ & $2.73\times10^{-3}$ & $7.88\times10^{-3}$ & m/s \\
\quad $\sigma_{r_x,t_f}$          & $2.23\times10^{-1}$ & $9.73\times10^{-1}$ & $3.37\times10^{-1}$ & $9.82\times10^{-1}$ & m \\
\quad $\sigma_{r_y,t_f}$          & $4.21\times10^{-1}$ & $8.48\times10^{-1}$ & $2.04\times10^{-1}$ & $7.17\times10^{-1}$ & m \\
\quad $\sigma_{v_x,t_f}$          & $1.19\times10^{-1}$ & $2.72\times10^{-1}$ & $1.36\times10^{-1}$ & $2.84\times10^{-1}$ & m/s \\
\quad $\sigma_{v_y,t_f}$          & $1.20\times10^{-1}$ & $2.94\times10^{-1}$ & $1.52\times10^{-1}$ & $2.99\times10^{-1}$ & m/s \\
\bottomrule
\end{tabular}
\end{threeparttable}
\end{table*}

\subsubsection{Gaussian Case}
\label{subsubsec:rktland_gaussian_results}

For the noise-free evaluation of the Gaussian case, the policy achieves mean final mass $\hat m_f=50{,}077.77~\mathrm{kg}$, corresponding to mean propellant consumption $\Delta \hat m=4{,}922.23~\mathrm{kg}$ and mean equivalent velocity increment $\Delta \hat v_{\mathrm{eq}}=407.45~\mathrm{m/s}$. Relative to the deterministic SCP baseline, this corresponds to an overhead of only $28.08~\mathrm{kg}$ in mean propellant and $2.44~\mathrm{m/s}$ in equivalent $\Delta v$, which is well below $1\%$.

\begin{figure*}[h!tb]
  \centering
  \includegraphics[width=0.98\textwidth]{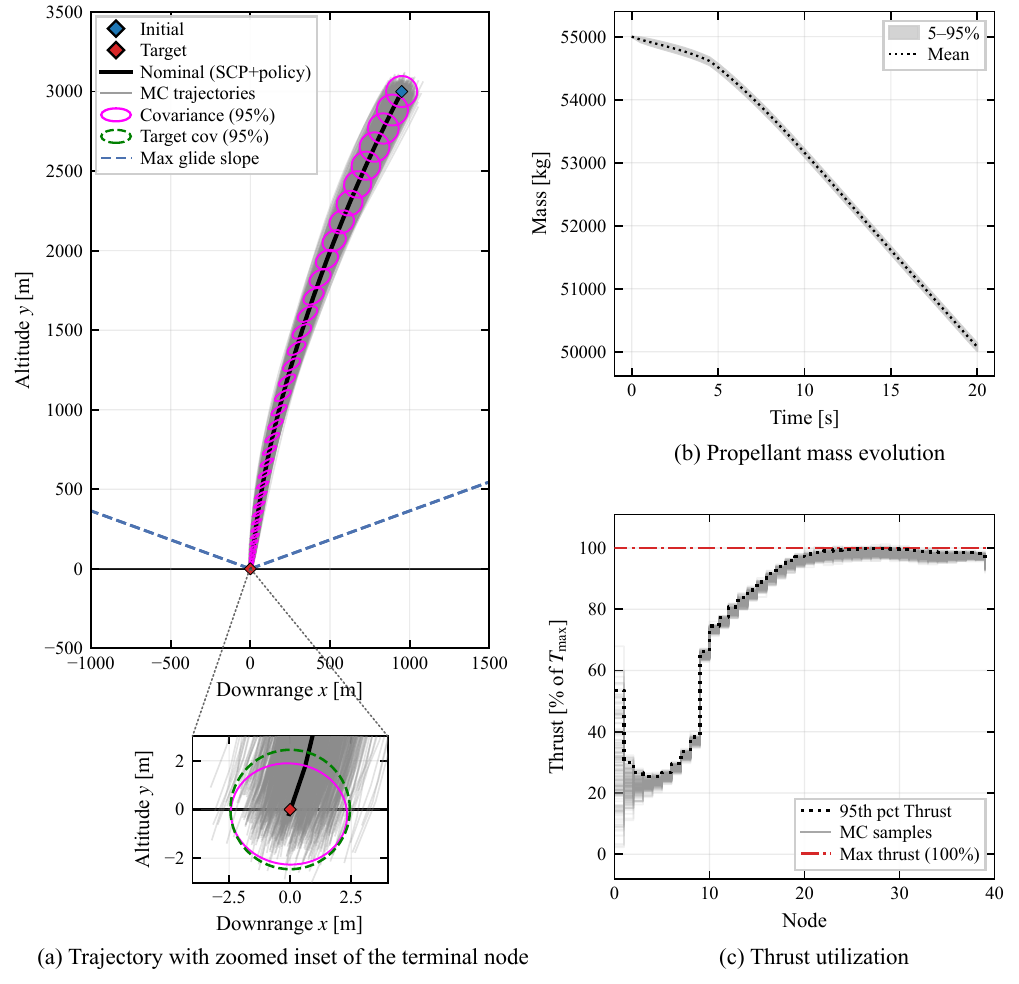}
  \caption{Closed-loop descent trajectory for the rocket landing case under Gaussian initial uncertainty with nominal process noise enabled. The trajectory bundle and covariance ellipses are shown with $4\times$ visual magnification, except for the terminal inset, which is plotted at true scale. Panels show: (a) trajectory bundle with nominal reference and terminal inset, (b) propellant mass evolution with mean and 5--95\% band, and (c) thrust utilization across the trajectory nodes.}
  \label{fig:rl_gaussian_traj}
\end{figure*}

\figref{fig:rl_gaussian_traj} shows the resulting closed-loop descent. The trajectory bundle remains tightly clustered around the nominal reference and contracts steadily toward touchdown. The thrust profile increases toward the cap during the terminal braking phase, but the enforced node-wise 95th-percentile profile remains below the available thrust across all nodes, indicating satisfaction of the implemented chance-constrained thrust specification.
Enabling nominal process noise during evaluation produces modest changes, indicating that for this descent profile and control schedule, process noise affects terminal accuracy far more strongly than fuel expenditure. The intended glide-slope geometry remains satisfied in both noise settings.

\begin{figure*}[h!tb]
  \centering
  \includegraphics[width=0.98\textwidth]{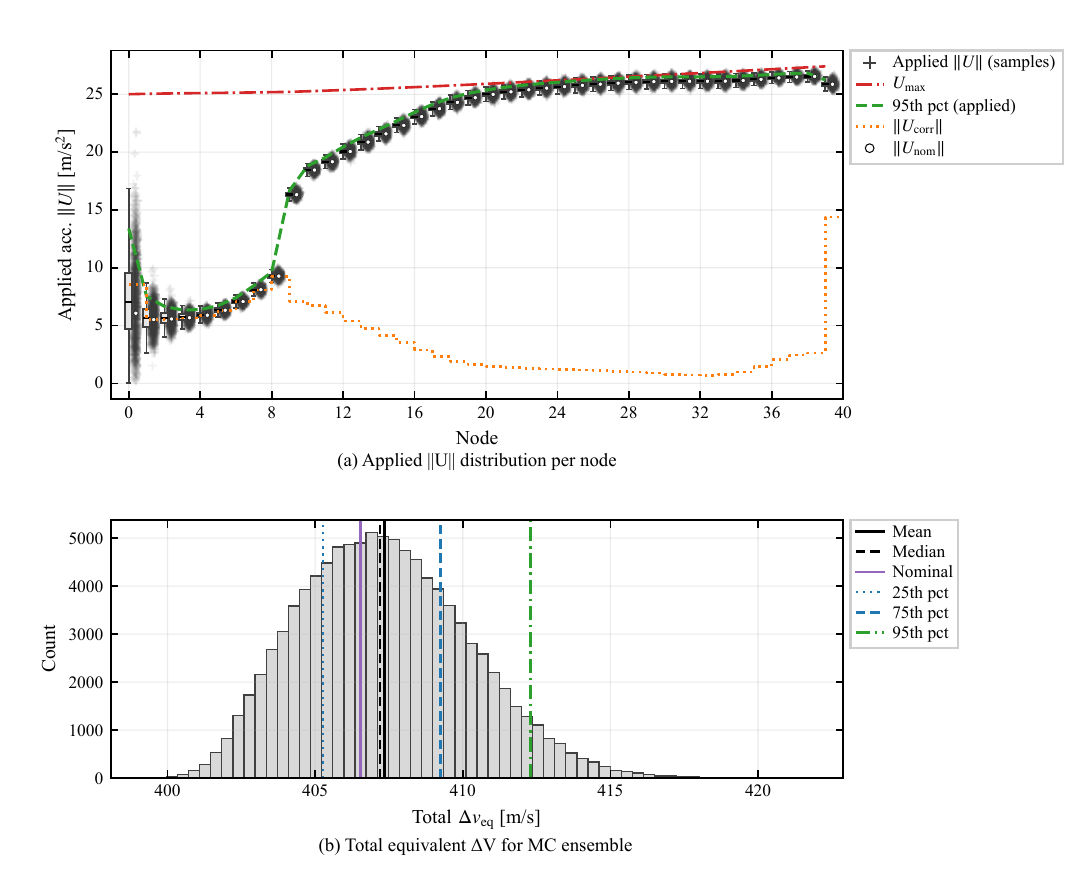}
  \caption{Applied control and equivalent-$\Delta v$ statistics for the rocket landing case under Gaussian initial uncertainty with nominal process noise enabled. Panels show: (a) node-wise distribution of applied acceleration magnitude $\norm{\vect{U}}$ across the evaluation ensemble and (b) histogram of total $\Delta v_{\mathrm{eq}}$ over the MC rollout.}
  \label{fig:rl_gaussian_ctrl}
\end{figure*}

\pagebreak

The control statistics in \figref{fig:rl_gaussian_ctrl} show that the applied $\norm{\vect{U}}$ distribution remains concentrated around a smooth profile, with the 95th-percentile staying below the thrust cap throughout the descent. The equivalent-$\Delta v$ distribution is also tightly concentrated, consistent with the narrow 5th--95th percentile spread reported in \tabref{tab:rl_results_comparison}.

\begin{figure*}[h!tb]
  \centering
  \includegraphics[width=0.98\textwidth]{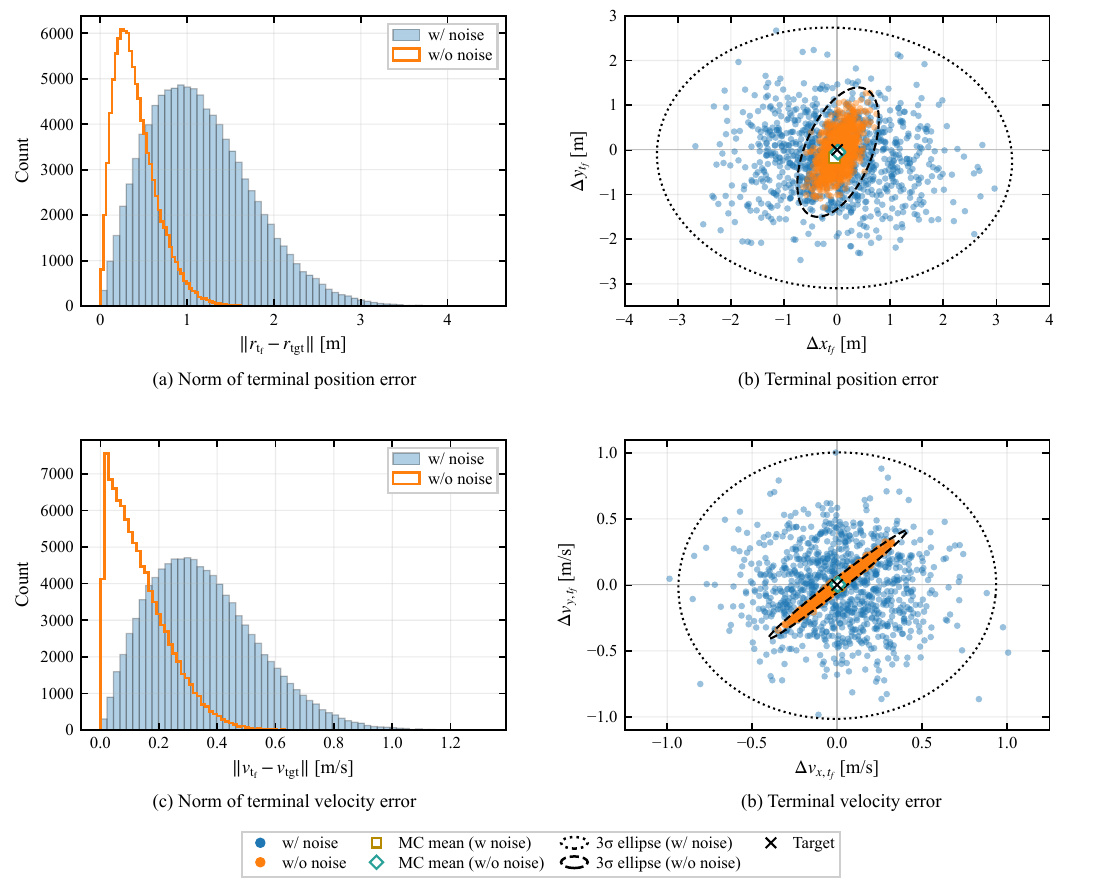}
  \caption{Terminal targeting error distributions for the rocket landing case under Gaussian initial uncertainty, shown both without and with nominal process noise. Panels show: (a) histogram of terminal position error norm, (b) terminal position-error components with mean and $3\sigma$ ellipses, (c) histogram of terminal velocity error norm, and (d) terminal velocity-error components with mean and $3\sigma$ ellipses.}
  \label{fig:rl_gaussian_disp}
\end{figure*}

\pagebreak

The terminal targeting statistics are reported in \tabref{tab:rl_results_comparison} and visualized in \figref{fig:rl_gaussian_disp}. 
In the noise-free setting, the mean touchdown bias is very small, 
$|\Delta \hat r_{x,t_f}|=1.76\times10^{-2}~\mathrm{m}$ and
$|\Delta \hat r_{y,t_f}|=5.22\times10^{-2}~\mathrm{m}$,
with terminal standard deviations
$\sigma_{r_x,t_f}=2.23\times10^{-1}~\mathrm{m}$,
$\sigma_{r_y,t_f}=4.21\times10^{-1}~\mathrm{m}$,
$\sigma_{v_x,t_f}=1.19\times10^{-1}~\mathrm{m/s}$, and
$\sigma_{v_y,t_f}=1.20\times10^{-1}~\mathrm{m/s}$.
These values indicate sub-meter positional spread and sub-$0.2~\mathrm{m/s}$ velocity spread at touchdown.
When nominal process noise is enabled at evaluation, the principal effect is instead on touchdown dispersion, which inflates to
$\sigma_{r_x,t_f}=9.73\times10^{-1}~\mathrm{m}$, 
$\sigma_{r_y,t_f}=8.48\times10^{-1}~\mathrm{m}$, 
$\sigma_{v_x,t_f}=2.72\times10^{-1}~\mathrm{m/s}$, and 
$\sigma_{v_y,t_f}=2.94\times10^{-1}~\mathrm{m/s}$.
Even in this noisy evaluation, the landing dispersion remains at sub-meter scale in position and below $0.3~\mathrm{m/s}$ in velocity. The covariance-feasibility indicator remains satisfied in both evaluations, so this result should be interpreted as a feasible but broadened terminal distribution.
The Gaussian landing case therefore shows the same qualitative behavior as the Earth--Mars process-noise study, i.e., disturbance mismatch broadens the terminal distribution far more than it affects aggregate fuel use.

\FloatBarrier

\subsubsection{Uniform Case}
\label{subsubsec:rktland_uniform_results}

\begin{figure*}[h!tb]
  \centering
  \includegraphics[width=0.98\textwidth]{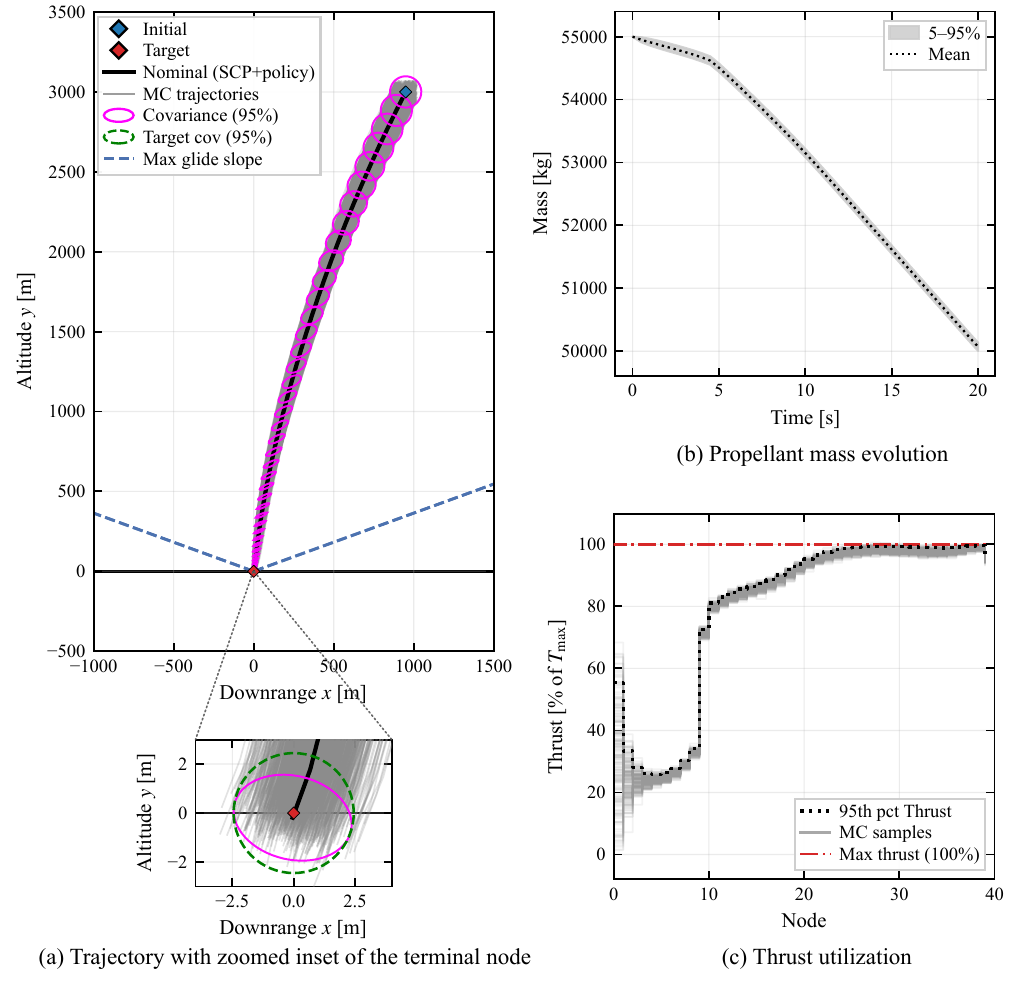}
  \caption{Closed-loop descent trajectory for the rocket landing case under bounded uniform initial uncertainty with nominal process noise enabled. The trajectory bundle and covariance ellipses are shown with $4\times$ visual magnification, except for the terminal inset, which is plotted at true scale. Panels show: (a) trajectory bundle with nominal reference, covariance ellipses, and glide-slope corridor, (b) propellant mass evolution with mean and 5--95\% band, and (c) thrust utilization across the trajectory nodes.}
  \label{fig:rl_uniform_traj}
\end{figure*}

To test distributional robustness beyond Gaussian assumptions, the policy is also trained and evaluated under bounded uniform initial-state uncertainty with matched first and second moments. \figref{fig:rl_uniform_traj} shows the resulting closed-loop descent trajectories under nominal process noise, where the same general behavior is observed as in the Gaussian case, with a tightly clustered trajectory bundle and a thrust profile that remains below the cap across the descent.

Under the noise-free evaluation, the uniform-trained policy exhibits slightly higher propellant usage than the Gaussian-trained policy, with $\Delta \hat m=4{,}928.24~\mathrm{kg}$ versus $4{,}922.23~\mathrm{kg}$, an increase of only $6.01~\mathrm{kg}$ ($0.122\%$). The corresponding mean equivalent velocity rises from $407.45$ to $407.97~\mathrm{m/s}$, and the 95th-percentile value from $412.84$ to $413.22~\mathrm{m/s}$. 
Enabling process noise at evaluation shows negligible change in mean fuel 
\pagebreak
usage relative to the corresponding noise-free result, with $\Delta \hat m=4{,}928.34~\mathrm{kg}$ and $\Delta \hat v_{\mathrm{eq}}=407.98~\mathrm{m/s}$. Relative to the noisy Gaussian evaluation, the uniform policy retains only a very small fuel penalty, with $\Delta v_{\mathrm{eq},95}=413.45~\mathrm{m/s}$ versus $412.85~\mathrm{m/s}$. The resulting node-wise 95th-percentile thrust remains below the cap throughout the descent, and the mean glide-slope corridor remains satisfied.

\begin{figure*}[h!tb]
  \centering
  \includegraphics[width=0.98\textwidth]{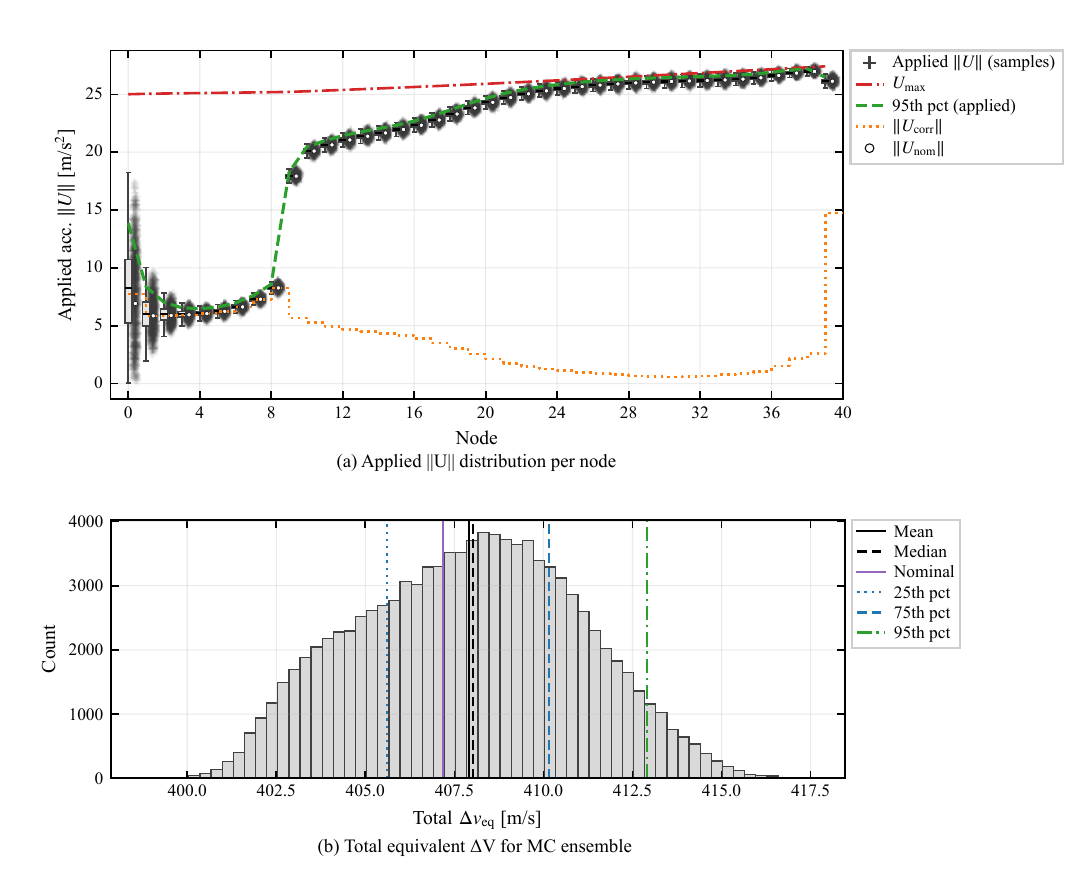}
  \caption{Applied control and equivalent-$\Delta v$ statistics for the rocket landing case under bounded uniform initial uncertainty with nominal process noise enabled. Panels show: (a) node-wise distribution of applied acceleration magnitude $\norm{\vect{U}}$ across the evaluation ensemble and (b) histogram of total $\Delta v_{\mathrm{eq}}$ over the MC rollout.}
  \label{fig:rl_uniform_ctrl}
\end{figure*}

\begin{figure*}[h!tb]
  \centering
  \includegraphics[width=0.98\textwidth]{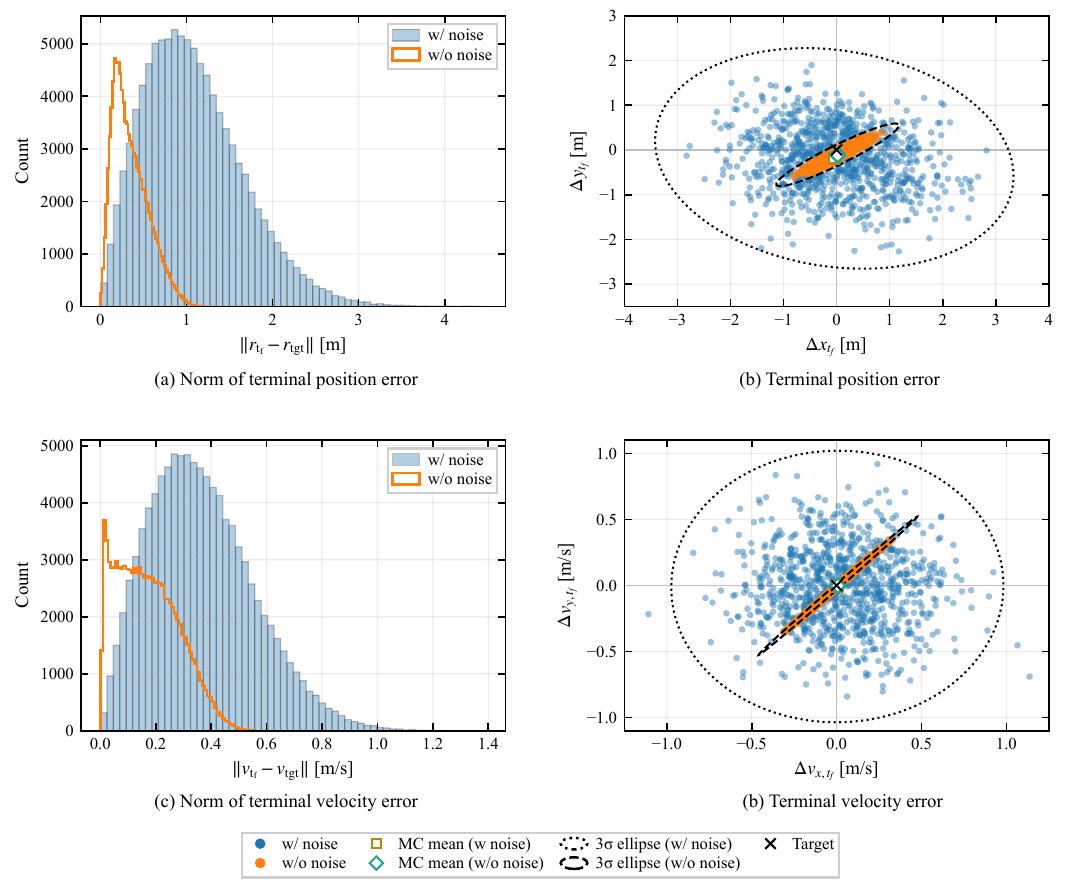}
  \caption{Terminal targeting error distributions for the rocket landing case under bounded uniform initial uncertainty, shown both without and with nominal process noise. Panels show: (a) histogram of terminal position error norm, (b) terminal position-error components with mean and $3\sigma$ ellipses, (c) histogram of terminal velocity error norm, and (d) terminal velocity-error components with mean and $3\sigma$ ellipses.}
  \label{fig:rl_uniform_disp}
\end{figure*}

The main qualitative difference between the rocket landing policies appears in the geometry of the terminal dispersion. In the noise-free setting, the uniform-trained policy yields
$\sigma_{r_x,t_f}=3.37\times10^{-1}~\mathrm{m}$, 
$\sigma_{r_y,t_f}=2.04\times10^{-1}~\mathrm{m}$, 
$\sigma_{v_x,t_f}=1.36\times10^{-1}~\mathrm{m/s}$, and 
$\sigma_{v_y,t_f}=1.52\times10^{-1}~\mathrm{m/s}$,
which differ from the Gaussian case not primarily in overall scale but in component-wise allocation. This produces a visibly different touchdown-ellipse geometry in \figref{fig:rl_uniform_disp}, indicating that the policy redistributes corrective authority differently under bounded non-Gaussian uncertainty even when fuel usage remains nearly unchanged. With nominal process noise enabled, the terminal spread inflates to
$\sigma_{r_x,t_f}=9.82\times10^{-1}~\mathrm{m}$, 
$\sigma_{r_y,t_f}=7.17\times10^{-1}~\mathrm{m}$, 
$\sigma_{v_x,t_f}=2.84\times10^{-1}~\mathrm{m/s}$, and 
$\sigma_{v_y,t_f}=2.99\times10^{-1}~\mathrm{m/s}$.
As in the Gaussian case, the dominant effect of process noise is therefore on touchdown dispersion rather than fuel consumption. Compared with the noisy Gaussian evaluation, the overall dispersion scale is similar, but modest anisotropy differences remain, which is consistent with the different geometry of the bounded initial-state model.

Taken together, the rocket-landing results show that RCCRL remains effective in a short-horizon continuous-thrust descent problem with drag, mass depletion, and glide-slope constraints, while preserving the same nominal-backbone, affine-correction, and sample-based chance-constraint structure used in the Earth--Mars transfer of Case Study~\Romannum{1}.
\section{Conclusion}
\label{sec:conclusion}

This paper presented a robust chance-constrained reinforcement learning (RCCRL) framework for spacecraft trajectory optimization in which a deterministic nominal trajectory is computed offline and then robustified through a learned affine closed-loop correction law. Rather than encoding probabilistic feasibility through problem-specific analytical reformulations, RCCRL evaluates chance constraints empirically through rollout ensembles and upper-tail quantile statistics, while terminal dispersion is regulated through covariance-based feasibility conditions applied to ensemble statistics. In this way, the framework preserves the structure and interpretability of a nominal trajectory backbone while extending it to a sample-based robust closed-loop setting.

The two case studies show that RCCRL is best understood as a reusable robustification architecture rather than as a controller tailored to a single mission class. Across both problems, what remains invariant is the same division of labor between deterministic nominal design and learned robust closed-loop augmentation, despite substantial differences in dynamical structure, control parameterization, horizon length, and constraint geometry. In the Earth--Mars transfer, the framework was benchmarked directly against a recent robust chance-constrained trajectory-optimization reference under Gaussian uncertainty and was then further examined under bounded uniform uncertainty and under process disturbances not seen during training. In the atmospheric pinpoint landing problem, the same methodological scaffold remained effective in a materially different regime involving short-horizon continuous-thrust descent, drag, mass depletion, and glide-slope constraints. Taken together, these results indicate that what transfers across problems is a coherent design logic built around deterministic nominal-reference generation, structured affine closed-loop augmentation, and sample-based evaluation of probabilistic feasibility.

This interpretation also clarifies the meaning of the distribution-agnostic RCCRL formulation and its role relative to classical robust optimization. Here, distribution-agnostic means that the control-design and evaluation methodology does not require uncertainty to belong to a closed-form family admitting a bespoke analytical reformulation, provided it can be sampled in simulation. It does not imply identical performance across uncertainty models. RCCRL therefore contributes a structured simulation-based methodology for robust trajectory optimization in regimes where nonlinear rollout evaluation, sampled uncertainty handling, and adaptive closed-loop correction are more natural than closed-form probabilistic reformulation. Its probabilistic guarantees are empirical rather than analytic, and their credibility therefore depends on simulation fidelity, uncertainty-model realism, and the accuracy of finite-sample quantile estimation. Within that trade space, however, the results of this paper support RCCRL as a practically useful and conceptually coherent methodology for robust spacecraft trajectory optimization.

Future work should extend the framework to incorporate higher-fidelity dynamics, richer disturbance and actuation models, and more explicit treatment of onboard state-estimation uncertainty, while also examining the trade-offs among nominal optimality, feedback authority, and probabilistic robustness across broader classes of spacecraft trajectory-design and optimization problems.

\section*{Acknowledgements}
This research was supported in part by the University of Auckland Doctoral Scholarship (Faculty of Engineering).
The authors wish to thank Jack Yarndley for helpful discussions on the development of the sequential convex programming setup.
The authors also wish to acknowledge the Centre for eResearch at the University of Auckland for their help in facilitating this research (URL: http://www.eresearch.auckland.ac.nz).


\end{document}